\documentclass%[draft]
{siamltex}
\usepackage{amsfonts}
\usepackage{amsmath}
\usepackage{amssymb}
\usepackage{array}
\usepackage{stmaryrd}
\usepackage{graphicx}
\usepackage{color}
\usepackage{tikz}
\usepackage{enumitem}

\newtheorem{remark}[theorem]{\it Remark}

\numberwithin{equation}{section}

\newcommand{\pt}{\partial}

\newcommand {\beq} {\begin{equation}}
\newcommand {\eeq} {\end{equation}}
\renewcommand{\sim}{\simeq}

\newcommand{\LL}{{\mathcal L}}
\renewcommand{\AA}{{\mathcal A}}

\newcommand{\RR}{{\mathcal R}}
\newcommand{\R}{\mathbb{R}}

\definecolor{pass}{rgb}{0,0,0.7}
\definecolor{green}{rgb}{0,0.6,0}
\definecolor{blue}{rgb}{0,0,0.8}
\begin{document}

\title{Error analysis for a fractional-derivative parabolic problem
on quasi-graded meshes\\
using barrier functions%
\thanks{The first author acknowledges support from Science Foundation Ireland Grant SFI/12/IA/1683.}%
}

\author{Natalia Kopteva%
\thanks{Department of Mathematics and Statistics, University of Limerick, Limerick, Ireland ({\tt natalia.kopteva@ul.ie}).}%
  \and Xiangyun Meng%
 \thanks{%
 Department of Mathematics, Beijing Jiaotong University, Beijing 100044, China ({\tt xymeng1@bjtu.edu.cn}).}%
 }

\date{}
\maketitle
\begin{abstract}
An initial-boundary value problem with a Caputo time derivative of fractional order $\alpha\in(0,1)$ is considered,
 solutions of which typically exhibit a singular behaviour at an initial time.
For  this problem, we
give a simple and general  numerical-stability analysis using barrier functions,
which yields sharp pointwise-in-time error bounds
on quasi-graded temporal meshes with arbitrary degree of grading.
L1-type and Alikhanov-type discretization in time are considered.
In particular, those results imply that milder (compared to the optimal) grading yields optimal convergence rates in positive time.
Semi-discretizations in time and full discretizations are addressed.
The theoretical findings are illustrated by numerical experiments.
\end{abstract}

\section{Introduction}
In this paper we give a simple and general numerical-stability analysis %of certain discretizations
for an initial-boundary value problem with a Caputo time derivative of fractional order $\alpha\in(0,1)$.
\smallskip
\begin{itemize}[leftmargin=0.7cm]
\item
The subtle and sharp stability property~\eqref{main_stab}, that we obtain, easily yields sharp pointwise-in-time error bounds for quasi-graded termporal meshes with
arbitrary degree of grading.
We are not aware of any such general results in the literature.

\item
In particular, our error bounds accurately
predict that milder (compared to the optimal) grading %that %coefficient threshold  for %
yields
optimal convergence rates in positive time.
This finding is new, and of practical importance.%(which has not been previously addressed).

\item
The simplicity of our approach is due to the usage of %clever/
versatile  barrier functions, which can
be used in the analysis of any discrete fractional-derivative operator that satisfies the discrete
maximum principle (or, more generally, is associated with an inverse-monotone matrix).

\item
Here  this approach is employed in the error analysis of the
L1 and Alikhanov L2-1${}_\sigma$  fractional-derivative operators, while in \cite{NK_L2} it is used in the analysis
of an L2-type discretization of order $3-\alpha$ in time.
In  \cite{NK_semil_L1} this methodology is generalized for semilinear fractional parabolic equations.

\end{itemize}
\smallskip

The Caputo fractional derivative in time, denoted here by $D_t^\alpha$, is  defined \cite{Diet10} by
\begin{equation}\label{CaputoEquiv}
D_t^{\alpha} u(\cdot,t) :=  \frac1{\Gamma(1-\alpha)} \int_{0}^t(t-s)^{-\alpha}\, \pt_s u(\cdot, s)\, ds
    \qquad\text{for }\ 0<t \le T,
\end{equation}
where $\Gamma(\cdot)$ is the Gamma function, and $\pt_s$ denotes the partial derivative in $s$.

Our main stability result  is that given an inverse-monotone fractional-derivative operator $\delta_t^\alpha$,
associated with a temporal mesh $\{t_j\}_{j=0}^M$ on $[0,T]$ with $\tau:=t_1$, and
 $\gamma\in \R$,   under certain conditions on the mesh, the following is true for $\{V^j\}_{j=0}^M$:
\beq\label{main_stab}
\left.\begin{array}{c}
|\delta_t^\alpha V^j|\lesssim (\tau/ t_j)^{\gamma+1}
\\[0.2cm]
\forall j\ge1,\;\;\; V^0=0
\end{array}\right\}
\;\;\Rightarrow\;\;
|V^j|\lesssim {\mathcal V}^j:=
\tau t_j^{\alpha-1}\left\{\begin{array}{ll}
1&\mbox{if~}\gamma>0\\
1+\ln(t_j/\tau)&\mbox{if~}\gamma=0\\
(\tau/t_j)^\gamma
%\tau^\gamb t_j^{\alpha-\gamb}
&\mbox{if~}\gamma<0
\end{array}\right.
\;\;\forall j\ge 1.
\eeq
This result is sharp in the sense that it
is consistent with the analogous property for the continuous Caputo operator $D_t^\alpha$; see Remark~\ref{rem_stab_sharp}.
The immediate usefulness of this property is due to the fact that truncation errors in time are typically bounded by negative powers of $t_j$.

It should be noted that while the explicit inverse of $D_t^\alpha$ is easily available, the proof of \eqref{main_stab} for any discrete operator is quite non-trivial.
As an alternative, discrete Gr\"{o}nwall inequalities were recently employed in the error analysis of L1- and Alikhanov-type schemes \cite{sinum18_liao_et_al,arxiv18_liao_et_al,sinum19_liao_et_al}.
However, this approach involves intricate evaluations and, furthermore, yields less sharp error bounds (see Remarks~\ref{rem_positive_time} and~\ref{rem_positive_time_A_par} for a more detailed discussion).
Our approach is entirely different and is substantially more concise as we obtain \eqref{main_stab} using clever barrier functions,
%(building on the ideas from \cite[Appendix~A]{NK_MC_L1}),
while the numerical results indicate that our error bounds are sharp in the pointwise-in-time sense.

The following fractional-order parabolic problem is considered:
\beq\label{problem}
\begin{array}{l}
D_t^{\alpha}u+\LL u=f(x,t)\quad\mbox{for}\;\;(x,t)\in\Omega\times(0,T],\\[0.2cm]
u(x,t)=0\quad\mbox{for}\;\;(x,t)\in\pt\Omega\times(0,T],\qquad
u(x,0)=u_0(x)\quad\mbox{for}\;\;x\in\Omega.
\end{array}
\eeq
This problem is posed in a bounded Lipschitz domain  $\Omega\subset\R^d$ (where $d\in\{1,2,3\}$).
The spatial operator $\LL$ here is a linear second-order elliptic operator:
\beq\label{LL_def}
\LL u := \sum_{k=1}^d \Bigl\{-\pt_{x_k}\!(a_k(x,t)\,\pt_{x_k}\!u) + b_k(x,t)\, \pt_{x_k}\!u \Bigr\}+c(x,t)\, u,
\eeq
with sufficiently smooth coefficients $\{a_k\}$, $\{b_k\}$ and $c$ in $C(\bar\Omega)$, for which we assume that $a_k>0$ in $\bar\Omega$,
and also either $c\ge 0$ or $c-\frac12\sum_{k=1}^d\pt_{x_k}\!b_k\ge 0$.

The first part of the paper is devoted to
L1-type schemes for problem \eqref{problem}, which employ the discetization of  $D^\alpha_tu$ defined, for $m=1,\ldots,M$, by
\vspace{-0.1cm}
\beq\label{delta_def}
\delta_t^{\alpha} U^m :=  \frac1{\Gamma(1-\alpha)} \sum_{j=1}^m \delta_t U^j\!\int_{t_{j-1}}^{t_j}\!\!(t_m-s)^{-\alpha}\, ds,
\qquad
\delta_t U^j:=\frac{U^j-U^{j-1}}{t_j-t_{j-1}},\vspace{-0.1cm}%
\eeq
when associated with the temporal mesh $0=t_0<t_1<\ldots <t_M=T$ on $[0,T]$.
The generality of our approach is demonstrated in the second part of the paper by extending the stability and error analysis to higher-order
Alikhanov-type  schemes \cite{alikh}.

Similarly to \cite{Brunner_MC85,Brunner_book,ChenMS_JSC,NK_MC_L1,sinum18_liao_et_al,arxiv18_liao_et_al,mclean_etal_NM07,Mustapha_etal_SINUM14,stynes_etal_sinum17},
%\cite{Mustapha_etal_IMANUM12,Mustapha_etal_SINUM13}},
our main interest will be in graded temporal meshes as they offer an efficient way of computing reliable numerical approximations of solutions singular at $t=0$,
which is typical for \eqref{problem}.
In particular, \cite{ChenMS_JSC,NK_MC_L1,Mustapha_etal_SINUM14,stynes_etal_sinum17}
 give global-in-time error bounds on graded meshes for problems of type \eqref{problem}
 for the L1 method \cite{stynes_etal_sinum17,NK_MC_L1}, the Alikhanov method \cite{ChenMS_JSC},
 and a high-order Petrov-Galerkin method in time \cite{Mustapha_etal_SINUM14}.
There is also a lot of interest in the literature  in optimal error bounds in positive time on uniform meshes;
see, e.g., \cite{gracia_etal_cmame,laz_L1,laz_review,NK_MC_L1}.
\smallskip

\begin{itemize}[leftmargin=0.7cm]
\item
By contrast,
here, as well as in the related paper \cite{NK_L2},  pointwise-in-time error bounds will be obtained, while
an arbitrary degree of mesh grading (with uniform meshes included as a particular case) is allowed.
%In particular, our results imply that milder (compared to the optimal) grading yields optimal convergence rates in positive time.
%and~\ref{rem_global_time}.}

\item
More general temporal meshes, which may be viewed as obtained by
adding new nodes in an arbitrary manner to a quasi-graded mesh, are also considered; see Section \ref{subs_gen_meshes} and Theorems \ref{lem_simplest_star} and \ref{theo_semidiscr}.

\item
For both considered discretizations,
when the optimal grading parameter $r=(p-\alpha)/\alpha$ is used, we recover the optimal global
convergence rates of $p-\alpha$, where $p=2$ for the L1 scheme and $p=3$ for the Alikhanov scheme,
as particular cases of our more general error bounds; see Remarks~%\ref{rem_positive_time},
\ref{rem_global_time} %\ref{rem_positive_time_A}
and~\ref{rem_global_time_A}.

\item
Another straightforward particular case of our error bounds indicates that
the optimal convergence rates  $p-\alpha$ in positive time $t\gtrsim 1$ are attained using much milder  grading
with $r>p-\alpha$;
%, while in the special case of the Alikhanov method applied to the parabolic problem
%a much milder grading $r=2$ yields the optimal convergence order $2$
see Remarks~\ref{rem_positive_time} and \ref{rem_positive_time_A}.
The accuracy of these threshold values is demonstrated by the  numerical results in \S\S\ref{ssec_num_L1}--\ref{ssec_num_A}.
%and \ref{rem_positive_time_A_par}).

\item
When dealing with the fractional parabolic case, for L1-type schemes we follow \cite{NK_MC_L1}, while for Alikhanov-type schemes, our approach substantially differs from \cite{alikh,ChenMS_JSC}
(as we aim at pointwise-in-time error bounds),
so we build on some ideas from \cite{NK_L2}, which  may be of independent interest.

\item
In the latter case, a much milder grading with $r=2$ (compared to the optimal $r=2/\alpha$) yields the optimal convergence order $2$; see Remark~\ref{rem_positive_time_A_par}.

\end{itemize}
\smallskip

Throughout the paper, it is assumed that there exists a unique solution of this problem %in $C(\bar\Omega\times[0,T])$
such that $\|\pt_t^l u(\cdot,t)\|_{L_2(\Omega)}\lesssim 1+t^{\alpha-l}$ for $l\le 3$.
%(the notation $\lesssim$ is rigourously defined in the final paragraph of this section).
This is a realistic assumption, satisfied by typical solutions of problem \eqref{problem}
(see, e.g.,   \cite{sakamoto}, %for smooth domains, \cite[\S2.2 and \S3.4]{laz_semidiscr} for polygonal/polyhedral domains,
\cite[\S2]{stynes_etal_sinum17}, \cite[\S6]{NK_MC_L1}),
in contrast to stronger assumptions of type $\|\pt^l u(\cdot,t)\|_{L_2(\Omega)}\lesssim 1$ frequently made in the literature
(see, e.g., references in \cite[Table~1.1]{laz_2fully_16}).
Indeed,
 \cite[Theorem~2.1]{stynes_too_much_reg} shows
that
if a solution $u$ of \eqref{problem} is less singular than we assume, %(in the sense that $|\pt_t^l u(\cdot,t)|\lesssim 1+t^{\gamma-l}$ for $l=0,1,2$ with any $\gamma>\alpha$),
then the initial condition $u_0$ is uniquely defined by the other data of the problem, which is clearly too restrictive.
At the same time, our results can be easily applied to the case of $u$ having no
singularities or exhibiting a somewhat different singular behaviour at $t=0$.
\smallskip

\begin{remark}\label{rem_stab_sharp}%[Sharpness of \eqref{main_stab}]
The stability result \eqref{main_stab} is sharp in the sense that it
is consistent with the analogous property for the continuous Caputo operator $D_t^\alpha$. Indeed, a calculation shows that
if $v(0)=0$ and $D_t^\alpha v(t)=F(t):=\min\{1, (\tau/t)^{\gamma+1}\}$ for $t>0$, then the explicit representation
$v(t)=J^\alpha_t F(t)=\{\Gamma(\alpha)\}^{-1}\!\! \int_{0}^t(t-s)^{\alpha-1} F( s)\, ds$
%as $J^{1-\alpha}v(t):=\{\Gamma(1-\alpha)\}^{-1}\!\! \int_{0}^t(t-s)^{-\alpha} v( s)\, ds$
 yields $v(t)\simeq {\mathcal V}(t)$ for $t\ge \tau$, where ${\mathcal V}(t)$ is a continuous version of ${\mathcal V}^j$ from \eqref{main_stab}.
\end{remark}
\smallskip

{\it Outline.} \S\ref{sec_L1_stab} is devoted to the proof of the stability result~\eqref{main_stab}
for the L1 discrete fractional-derivative operator. This result is then employed in \S\ref{sec_L1_error} to obtain
pointwise-in-time error bounds for
L1-type discretizations of the initial-value problem in
 \S\ref{ssec_L1_prdgm}, as well as semi-discretizations and full discretizations of the initial-boundary-value problems in \S\ref{ssec_L1_semi} and \S\ref{ssec_L1_full}.
 The above error analysis is extended to the Alikhanov-type discretizations in
\S\ref{sec_A_stab}.
Finally, our theoretical findings are illustrated by numerical experiments in \S\ref{sec_Num}.
\smallskip

{\it Notation.}
We write
 $a\sim b$ when $a \lesssim b$ and $a \gtrsim b$, and
$a \lesssim b$ when $a \le Cb$ with a generic constant $C$ depending on $\Omega$, $T$, $u_0$ and
%$C_f$ (and possibly $\bar C_f$),
$f$,
but not %on other essential quantities.
%
 %In particular,
% $C$  does not depend
 on the total numbers of degrees of freedom in space or time.
  Also, for %$\mathcal{D}\subset\bar\Omega $,
  $1 \le p \le \infty$, and $k \ge 0$,
  we shall use the standard norms
% $\|\cdot\|_{L_p(\mathcal{D})}$
%  and $\|\cdot\|_{W_p^k(\mathcal{D})}$ respectively
  in the spaces $L_p(\Omega)$ and the related Sobolev spaces $W_p^k(\Omega)$,
  while  $H^1_0(\Omega)$ is the standard space of functions in $W_2^1(\Omega)$ vanishing on $\pt\Omega$.

\section{Stability properties of the L1 discrete fractional-derivative operator}\label{sec_L1_stab}

\subsection{Quasi-graded temporal meshes. Main stability result}

Throughout the paper, we shall frequently assume that the temporal mesh is quasi-graded
in the sense that, with some $r\ge 1$,
\beq\label{t_grid_gen}
\tau := t_1\simeq M^{-r},\quad
t_j \sim \tau j^r,
\quad %t_j = \tau j^r,\quad
\tau_j:=t_j-t_{j-1}%\simeq M^{-1}\,t_j^{1-1/r}
\simeq\tau^{1/r}t_j^{1-1/r}%{\color{red}\simeq t_j/j}
\quad %t_{j+1}\lesssim  t_{j}
\forall\,j=1,\ldots,M.
\eeq

For example,  the standard graded temporal mesh
$\{t_j=T(j/M)^r\}_{j=0}^M$ with some $r\ge 1$ (while $r=1$ generates a uniform mesh)
 satisfies \eqref{t_grid_gen},
%For this mesh, a calculation shows that
%\beq\label{t_grid}
%\tau := t_1\simeq M^{-r},\quad t_j = \tau j^r,\quad \tau_j:=t_j-t_{j-1}%\simeq M^{-1}\,t_j^{1-1/r}
%\simeq\tau^{1/r}t_j^{1-1/r}
%\quad
%\mbox{for\;\;}j=1,\ldots,M.
%\eeq
in view of $\tau_j\simeq M^{-1}\,t_{j-1}^{1-1/r}$ %for $j=1$,
and $t_j\le 2^{r}  t_{j-1}$
for $j\ge 2$.

Furthermore, our results
also apply to more general %temporal
meshes that may be viewed as obtained
by adding new nodes to any mesh of
 type \eqref{t_grid_gen}; see Section~\ref{subs_gen_meshes}.

%\subsection{Main stability result}

The key  in our error analysis for L1-type discretizations is the following stability property.
\smallskip

\begin{theorem}[Stability]\label{theo_barrier}
(i) Let the temporal mesh satisfy
 \eqref{t_grid_gen}
with $1\le r\le (2-\alpha)/\alpha$.
Given  $\gamma\in \R$ and $\{V^j\}_{j=0}^M$,
% with $V^0=0$ one has
the stability property \eqref{main_stab} holds true.
%for $j=1,\ldots, M$.\\
\\[0.1cm]
(ii)
If $\gamma\le \alpha-1$, then \eqref{main_stab} holds true on an arbitrary temporal mesh $\{t_j\}_{j=0}^M$.
\\[0.1cm]
(iii) The above results remain valid if  $|\delta_t^\alpha V^j|\lesssim (\tau/ t_j)^{\gamma+1}$ in \eqref{main_stab} is replaced by $\delta_t^\alpha |V^j|\lesssim (\tau/ t_j)^{\gamma+1}$.%
\end{theorem}
\smallskip

\begin{proof}
%{\it Proof of Theorem~\ref{theo_barrier}.}
(i) It suffices
 to prove part (i) only for $\gamma\ge \alpha-1$ (as the result of part (ii) applies to the case
$\gamma\le \alpha-1$).
The proof is presented in Sections~\ref{ssec_B_alpha} and~\ref{ssec_2_2}, where a few cases are considered separately.
\smallskip

(ii)
This result is easily obtained from \cite[Lemma 2.1(i)]{NK_MC_L1}.
The latter implies that
$|V^m|\lesssim\max_{j\le m}\bigl\{t_j^\alpha|\delta_t^\alpha V^j|\bigr\}$ on an arbitrary mesh.
The assumptions on $\{V^j\}$ yield
$t_j^\alpha|\delta_t^\alpha V^j|\lesssim t_j^\alpha (\tau/ t_j)^{\gamma+1}
=\tau^{\gamma+1}  t_j^{\alpha-\gamma-1}$, which, combined with $\gamma\le \alpha-1$,
implies
$t_j^\alpha|\delta_t^\alpha V^j|
\lesssim \tau^{\gamma+1}  t_m^{\alpha-\gamma-1}$ $\forall\, j\le m$.
The desired assertion $|V^m|\lesssim  \tau^{\gamma+1}  t_m^{\alpha-\gamma-1}$ follows.
\smallskip

(iii) Imitate the proof of  \cite[Lemma 2.1(ii)]{NK_MC_L1}.
To be more precise, let $W^0=0$ and $\delta_t^\alpha W^j= \max\bigl\{0,\,\delta_t^\alpha |V^j|\bigr\}\ge \delta_t^\alpha |V^j|$  $\forall\,j\ge1$.
Then $0\le |V^j|\le W^j$ $\forall\,j\ge1$
(as $\delta_t^\alpha$ is associated with an $M$-matrix), while  the results of parts (i) and (ii) apply to $\{W^j\}$.
\end{proof}
%\smallskip

\begin{remark}
To a degree, the proof of Theorem \ref{theo_barrier} builds on the %stability
analysis in \cite[Appendix A]{NK_MC_L1} for uniform grids% %in \cite{NK_MC_L1}.
%(in particular,
, but now we address considerably more general meshes.%
\end{remark}

\subsection{Proof of Theorem~\ref{theo_barrier}(i) for
 $\gamma \ge \alpha$}\label{ssec_B_alpha}
In this case $(\tau/ t_j)^{\gamma+1}\le (\tau/ t_j)^{\alpha+1}$,
so it suffices to consider only $\gamma=\alpha$.
For the latter case,
as the operator $\delta_t^\alpha$ is associated with an $M$-matrix, it suffices to prove the following lemma.
\smallskip

\begin{lemma}\label{lem_stability1}
Let the temporal mesh satisfy
 \eqref{t_grid_gen}
with $1\le r\le (2-\alpha)/\alpha$.
Then there exists a discrete barrier function $\{B^j\}_{j=0}^M$  such that
\beq\label{B_def}
B^0=0,\qquad
0\le B^j\lesssim t_j^{\alpha-1},\qquad
\delta_t^\alpha B^j\gtrsim \tau^\alpha t_j^{-\alpha-1}
\qquad\forall j\ge1.
\eeq
\end{lemma}\vspace{-0.3cm}
%\smallskip
\begin{proof}
Fix a sufficiently large number $2\le p\lesssim 1$, and then
set %$\beta:=1-\alpha$ and %$B(s):=\min\bigl\{(s/t_p)t_p^{-\beta}, s^{-\beta} \bigr\}$, and also $B^j:=B(t_j)$.
\beq\label{B_L1_def}
\beta:=1-\alpha,\qquad
B(s):=\min\bigl\{(s/t_p)t_p^{-\beta}, s^{-\beta} \bigr\},\qquad B^j:=B(t_j).
\eeq
Note that, when using the notation of type $\lesssim$, the dependence on $p$ will be shown explicitly.

For $j\le p$, a straightforward calculation shows that
$\delta_t^\alpha B^j= D^\alpha_t B(t_j)\sim t_j^\beta t_p^{-\beta-1}\sim p^{-r(2-\alpha)}t_j^\beta \tau^{-\beta-1} $,
where we also used $t_p \sim \tau p^r$ (in view of \eqref{t_grid_gen}). As $t_j\ge \tau$, we then get
$\delta_t^\alpha B^j \gtrsim p^{-r(2-\alpha)} \tau^{-\beta} t_j^{\beta-1}\ge p^{-r(2-\alpha)}\tau^\gamma t_j^{-\gamma-1}$ $\forall \gamma\ge \alpha-1$
including $\gamma=\alpha$.

Next, for $D_t^\alpha B(t)$ with $t>t_p$ one has
\begin{align*}
\Gamma(1-\alpha)\,D^\alpha_t B(t)=&\underbrace{\int_0^{t_p}\!\! t_p^{-\beta-1}(t-s)^{-\alpha}\,ds}_%
{{}\ge t_p^{-\beta} t^{-\alpha}}
-\underbrace{\beta\int_{t_p}^t\! s^{-\beta-1}(t-s)^{-\alpha}\,ds}_%
{{}=:  t^{-1}{ I} }\,.
%\\\ge&
%t_1^{-\beta} t^{-\alpha} - t^{-1}\beta\int_{\hat t_1}^1 \hat s^{-\beta-1}(1-\hat s)^{-\alpha}\,d\hat s
\end{align*}
Here, using $\hat s:=s/t$ and $\hat t_p:=t_p/t$, and noting that $\alpha+\beta=1$, one gets
$$
{I}=\beta\int_{\hat t_p}^1 \hat s^{-\beta-1}(1-\hat s)^{-\alpha}\,d\hat s
=\hat t_p^{-\beta}(1-\hat t_p)^{\beta}
\le \hat t_p^{-\beta}(1-\beta \hat t_p).
$$
Now, using $t^{-1}\hat t_p^{-\beta}=t_p^{-\beta} t^{-\alpha}$, one concludes %for $t>t_p$
that
\beq\label{app_B_eq}
\Gamma(1-\alpha)\,
D^\alpha_t B(t)
%\ge t_p^{-\beta} t^{-\alpha}\Bigl\{1-(1-\hat t_p)^{1-\alpha}\Bigr\}
\ge
%\{\Gamma(1-\alpha)\}^{-1}
t_p^{-\beta} t^{-\alpha}\,(\beta t_p/t)
=%\{\Gamma(1-\alpha)\}^{-1}
\beta t_p^{\alpha}t^{-\alpha-1}
\qquad\mbox{for}\;\;t>t_p.
%={\color{blue}\beta p^{r\alpha}\,(\tau^\alpha t^{-\alpha-1})}.
\eeq
So, to complete the proof,
 it remains to show that $|\delta^\alpha_t B^m-D^\alpha_t B(t_m)|\le \frac12 D^\alpha_t B(t_m)$ for any $m> p$.

For the latter,
note that   $\Gamma(1-\alpha)[\delta^\alpha_t B^m-D^\alpha_t B(t_m)]=\sum_{j=1}^m\mu^j$, where,
using the standard piecewise-linear interpolant $B^I$ of $B$,
\beq\label{mu_def}
\mu^j:= %\int_{t_{j-1}}^{t_j}\!\![\delta_t B^j-B'(s)]\underbrace{(t_m-s)^{-\alpha}}_{{}=g(s)}ds=
\int_{t_{j-1}}^{t_j}\!\!(B^I-B)'(s)\,(t_m-s)^{-\alpha}ds
=\alpha\int_{t_{j-1}}^{t_j}\!\! (B-B^I)(s)\,(t_m-s)^{-\alpha-1} ds.
\eeq
Clearly, $\mu^j=0$ for $j\le p$.
For $p+1\le j\le m-1$,
one gets $| B-B^I |\lesssim \tau_j^2 |B''(t_{j-1})|$.
%$\lesssim \tau_j^2 |B''(s)|\sim \tau_j^2 s^{-\beta-2}$ for $s\in(t_{j-1},t_j)$ (in view of $t_{j-1}\sim t_j$).
For $j=m$, we shall use a similar, but sharper bound $| B-B^I |\lesssim  \tau_j(t_m-s) |B''(t_{j-1})|$.
Combining these  yields $|B-B^I |\le \tau^2_j\min\bigl\{1,(t_m-s)/\tau_m\bigr\} |B''(t_{j-1})|$ for $j>p$, where
$|B''(t_{j-1})|\lesssim  |B''(s)|\sim  s^{-\beta-2}$ (in view of $t_{j-1}\sim t_j$).
Also noting that, in view of~\eqref{t_grid_gen}, $\tau_j\simeq\tau^{1/r}t_j^{1-1/r}\simeq\tau^{1/r}s^{1-1/r}$, we arrive at
$$
|\mu^j| %= \left|\int_{t_{j-1}}^{t_j}(B-B^I) f'ds \right|
%\lesssim \left|\int_{t_{j-1}}^{t_j}\tau_j^2 {\color{red}\min\bigl\{(1,(t_m-s)/\tau_m\bigr\}}B'' g'ds\right|
\lesssim \tau^{2/r}\!\int_{t_{j-1}}^{t_j}\!\!\!s^{-\beta-2/r}\,(t_m-s)^{-\alpha-1}\,\min\bigl\{1,(t_m-s)/\tau_m\bigr\}\,ds
\qquad\forall\,p>m.
$$
This immediately yields the bound
\beq\label{mu_sum}
\left|\delta^\alpha_t B^m-D^\alpha_t B(t_m)\right|\lesssim
 \tau^{2/r}\!\int_{t_{p}}^{t_m}\!\!\!s^{-\beta-2/r}\,(t_m-s)^{-\alpha-1}\,\min\bigl\{1,(t_m-s)/\tau_m\bigr\}\,ds.
\eeq
For the latter, using the substitution $s=t_m \hat s$ and the notation $\hat t_j:=t_j/t_m$, $\hat \tau_j:=\tau_j/t_m$, one gets
$$
\left|\delta^\alpha_t B^m-D^\alpha_t B(t_m)\right|\lesssim
\tau^{2/r}t_m^{-2/r-1}
\underbrace{\int_{\hat t_{p}}^{1}\!\hat s^{-\beta-2/r}(1-\hat s)^{-\alpha-1}\min\bigl\{1,(1-\hat s)/\hat\tau_m\bigr\}\,d\hat s}_{\lesssim\, \hat t_p^{\alpha-2/r}  +\hat\tau_m^{-\alpha}}.
$$
Here, when bounding the integral, it is convenient to separately consider the intervals
$(\hat t_p,\max\{\frac12, \hat t_p\})$, $(\max\{\frac12, \hat t_p\},1-\hat\tau_m)$ and $(1-\hat\tau_m,1)$,
where $\hat\tau_m\le \frac12$ if $p$ is sufficiently large (as, in view of~\eqref{t_grid_gen}, $\hat\tau_m\sim 1/m\le 1/p$).
On these three intervals, the integrand is respectively
$\lesssim s^{-\beta-2/r}$, $\lesssim (t_m-s)^{-\alpha-1}$ and $\lesssim (t_m-s)^{-\alpha}/\tau_m$, so the corresponding integrals are respectively
$\lesssim \hat t_p^{\alpha-2/r}$, $\lesssim \hat\tau_m^{-\alpha}$ and $\lesssim \hat\tau_m^{-\alpha}$.
Finally, note that $\hat\tau_m=\tau_m/t_m\simeq(\tau/t_m)^{1/r}$,
while, in view of $r\le (2-\alpha)/\alpha$, one has
$(\tau/t_m)^{(2-\alpha)/r}\lesssim (\tau/t_m)^{\alpha}$. Now, a calculation shows that
\begin{align}\notag
\left|\delta^\alpha_t B^m-D^\alpha_t B(t_m)\right|&\lesssim  \tau^{2/r}t_m^{-2/r-1}\Bigl[ (t_p/t_m)^{\alpha-2/r}  +(\tau/t_m)^{-\alpha/r}\Bigr]\\[0.3cm]\notag
&\lesssim  (\tau/t_p)^{2/r} t_p^\alpha t_m^{-\alpha-1}
+ \underbrace{t_m^{-1}(\tau/t_m)^{(2-\alpha)/r}}_{\lesssim \tau^\alpha t_m^{-\alpha-1}}\\
\label{mu_sum1}
&\lesssim
\Bigl[(\tau/t_p)^{2/r}+(\tau/t_p)^\alpha\Bigr]
\,t_p^\alpha t_m^{-\alpha-1}.
\end{align}
Combining this with \eqref{app_B_eq} and choosing $p$ sufficiently large yields the desired assertion
$\delta^\alpha_t B^m\gtrsim t_p^\alpha t_m^{-\alpha-1}$ $\forall\,m>p$, and hence $\forall\,m\ge1$.
\end{proof}
\smallskip

\begin{corollary}\label{cor_gen_mesh}
Lemma~\ref{lem_stability1} remains valid if
the temporal mesh is obtained by adding new  nodes to any mesh of type
 \eqref{t_grid_gen} under the condition that the first mesh interval remains unchanged.
\end{corollary}
\smallskip
\begin{proof}
Suppose the temporal mesh $\{t_k'\}$ is obtained by refining
the mesh $\{t_j\}_{j=0}^M$ of type \eqref{t_grid_gen}.
For $t_k'\le t_p$,
it is essential that $t_1'=t_1=\tau$,
so the desired result is obtained exactly as in the proof of Lemma~\ref{lem_stability1}.
For $t_k'> t_p$,
the desired result is obtained by combining
\eqref{app_B_eq}
with the bound $|D^\alpha_t (B^I-B)|\le \frac12 D^\alpha_t B$
at any $t_k'>t_p$, where $B^I$ denotes the piecewise-linear interpolant on the new mesh $\{t_k'\}$.
If $t_k'=t_m$ for some $m>p$, %the latter is obtained again using \eqref{mu_def}, in which %$B^I$ is the interpolant on the new mesh $\{t_k'\}$,
we again proceed exactly as in the proof of Lemma~\ref{lem_stability1},
as the same bounds on $B-B^I$ hold true (even though $B^I$ is now the interpolant on a finer mesh).
If $t_k'\in(t_m,t_{m+1})$ for $m\ge p$,
then on $(t_{m-1},t_k')$ one employs
$| B-B^I| \lesssim  \tau_m(t_k'-s) |B''(t_{m-1})|\lesssim  \tau_m^2\min\{1,\,(t_k'-s)/\tau_m\} |B''(t_{m-1})|$.
Hence, one gets a version of
\eqref{mu_sum} with $t_m$ replaced by $t_k'$, which (in view of $t_k'\simeq t_m$) leads to the desired version of \eqref{mu_sum1} at $t_k'$.
\end{proof}

\subsection{Proof of Theorem~\ref{theo_barrier}(i) for
 $\gamma < \alpha$}\label{ssec_2_2}

 We shall use the notation and some findings from the proof of Lemma~\ref{lem_stability1}.
 In particular,  $\beta=1-\alpha$, while $p\sim 1$ was chosen sufficiently large in the proof of Lemma~\ref{lem_stability1}.
When using the notation of type $\lesssim$, the dependence on $\gamma$ and $m$, but not on $p$, will be shown explicitly.

For $m\ge 0$ and $\gamma < \alpha$, set
\beq\label{p_m_c_m_notation}
p_m:=2^m p,
\qquad
B_m^j:=\min\bigl\{t_j t_{p_m}^{-\beta-1}, t_j^{-\beta} \bigr\},
%B_m(s):=\min\bigl\{s t_{p_m}^{-\beta-1}, s^{-\beta} \bigr\},
%\quad  B_m^j:=B_m(t_j),
\qquad
c_m:=2^{-m\gamma r}\;\;\Rightarrow\;\;
c_mt_{p_m}^\gamma\sim\tau^\gamma.
\eeq
Here the final observation follows from \eqref{t_grid_gen} (which yields $t_{p_m} \sim \tau p_m^r$).

Note that $B_0^j=B^j$, and, more generally, $B_m^j=B^j\big|_{p:=p_m}$, where $B^j$ is
from \eqref{B_L1_def}.
Conveniently, in the proof of Lemma~\ref{lem_stability1},
 the dependance on any sufficiently large $p$ was shown explicitly.
In particular, we recall that $\delta_t^\alpha B_m^j\ge 0$ for $j\ge 0$.
Furthermore,
\beq\label{B_m_j_bounds}
\delta_t^\alpha B_0^j\gtrsim \tau^\gamma t_j^{-\gamma-1}
\;\;\mbox{for}\;1\le j\le p_0,
\qquad
c_m(\delta_t^\alpha B_m^j)\gtrsim %c_m t_{p_m}^\gamma t_j^{-\gamma-1}\sim
\tau^\gamma t_j^{-\gamma-1}\;\;\mbox{for}\;p_m< j\le p_{m+1}.
\eeq
The first relation for $B_0^j=B^j$ can be found in the above-mentioned proof
for $\gamma\ge \alpha-1$ (but is, in fact, valid for any fixed $\gamma$ now that
the dependence on $p$ is inessential).
The second relation in \eqref{B_m_j_bounds} follows from the bound of type
\eqref{app_B_eq} also obtained there:
$\delta^\alpha_t B^j_m\gtrsim t_{p_m}^{\alpha}t_j^{-\alpha-1}$.
The latter, indeed, implies
$c_m(\delta_t^\alpha B_m^j)\gtrsim c_m t_{p_m}^\gamma t_j^{-\gamma-1}\sim \tau^\gamma t_j^{-\gamma-1}$ for $p_m< j\le p_{m+1}$
(also using the final bound in \eqref{p_m_c_m_notation}).

Now we are prepared to prove the following two lemmas, which are sufficient for establishing Theorem~\ref{theo_barrier}(i) for
 $\gamma \in(0, \alpha)$ and $\gamma\le 0$ respectively.
 \smallskip

% \begin{lemma}\label{lem_stability2}
%Given $0 < \gamma < \alpha$, then $1\le r < 2-\alpha$.
%For any $\{V^j\}_{j = 0}^M$ on graded mesh $\{t_j\}_{j = 0}^M$.
%If $V^0=0$ and $|\delta_t^\alpha V^j|\lesssim \tau^\gamma t_j^{-\gamma-1}$ for $j=1,\ldots,M$, then
%$|V^j|\lesssim t_j^{\alpha-1}$ for $j=1,\ldots, M$.
%\end{lemma}
 \begin{lemma}\label{lem_stability2}
Under the conditions of Theorem~\ref{theo_barrier}(i), suppose that $\gamma \in(0, \alpha)$. Then there exists a discrete barrier function $\{\bar B^j\}_{j=0}^M$  such that
$\bar B^0=0$, while
$0\le \bar B^j\lesssim t_j^{\alpha-1}$ and $\delta_t^\alpha \bar B^j\gtrsim \tau^\gamma t_j^{-\gamma-1}$
for $j\ge1$.
\end{lemma}
\smallskip
\begin{proof}
Using \eqref{p_m_c_m_notation}, let $\bar B^j:=\sum_{m=0}^\infty c_m B_m^j$.
Then $\delta_t^\alpha \bar B^j\gtrsim \tau^\gamma t_j^{-\gamma-1}$ $\forall j\ge 1$ follows from \eqref{B_m_j_bounds},
while
$\sum_{m=0}^\infty c_m= C_\gamma:=(1-2^{-\gamma r})^{-1}$,
so
$\bar B^j\le  C_\gamma t_j^{-\beta}=C_\gamma t_j^{\alpha-1}$, which   completes the proof.
\end{proof}
\smallskip
\begin{lemma}\label{lem_stability3}
Under the conditions of Theorem~\ref{theo_barrier}(i), suppose that $\gamma \in[\alpha-1, 0]$.
If $V^0=0$ and $|\delta_t^\alpha V^j|\lesssim \tau^\gamma t_j^{-\gamma-1}$ for $j=1,\ldots,n\le M$, then
$|V^n|\lesssim t_n^{\alpha-1}[1+\ln (t_n/\tau)]$ if $\gamma=0$,
and $|V^n|\lesssim t_n^{\alpha-1}(\tau/t_{n})^\gamma$ if $\gamma \le 0$.
\end{lemma}
\smallskip
\begin{proof}
Using \eqref{p_m_c_m_notation}, let $\bar B^j:=\sum_{m=0}^N c_m B_m^j$,
where $N=0$ if $n\le p$, and
$N := \lceil\log_2(n/p) -1\rceil$ otherwise,
so that $p_N<n\le p_{N+1}$.
Note also that $N\lesssim\ln n\simeq \ln (t_n/\tau)$ (as $t_n/\tau\sim n^r$ in view of \eqref{t_grid_gen}).
Then $\delta_t^\alpha \bar B^j\gtrsim \tau^\gamma t_j^{-\gamma-1}$ for $ 1\le j\le n$ follows from \eqref{B_m_j_bounds}.
Hence $|V^j|\lesssim \bar B^j$  $\forall j\le n$, in particular $|V^n|\lesssim \bar B^n$.

On the other hand,
$\bar B^n\le   t_n^{\alpha-1}\sum_{m=0}^N c_m$.
When $\gamma=0$, each $c_m=1$, so $\sum_{m=0}^N c_m=1+N\sim 1+\ln (t_n/\tau)$, so, indeed,
$\bar B^n\le t_n^{\alpha-1}[ 1+\ln (t_n/\tau)]$.
When $\gamma \in[\alpha-1, 0)$, we get
$\sum_{m=0}^N c_m=(c_{N+1}-1)/(c_1-1)$,
where
$c_{N+1}\sim (\tau/t_{p_{N+1}})^\gamma \sim(\tau/t_{n})^\gamma$,
while
$C_\gamma:=(c_1-1)^{-1}=(2^{|\gamma| r}-1)^{-1}$, so finally
$|V^n|\lesssim\bar B^n\lesssim C_\gamma t_n^{\alpha-1}(\tau/t_{n})^\gamma$.
\end{proof}

\subsection{More general temporal meshes}\label{subs_gen_meshes}
Our main stability result,  Theorem~\ref{theo_barrier}, remains valid
for more general temporal meshes that may be viewed as obtained
by adding new nodes to any mesh of
 type \eqref{t_grid_gen}
 under the condition that the first mesh interval remains unchanged.
 Indeed, an inspection of the proof of Corollary~\ref{cor_gen_mesh} reveals
that not only Lemma~\ref{lem_stability1}, but also the results of Section~\ref{ssec_2_2} are valid for the above temporal mesh.
Such more general meshes may be useful if the solution exhibits additional singularities away from the initial time.
\smallskip

\begin{lemma}
Theorem~\ref{theo_barrier} remains valid if the temporal mesh satisfies the following weaker version of \eqref{t_grid_gen}:
\beq\label{t_grid_gen_gen}
\tau := t_1\simeq M^{-r},\qquad
%t_j \sim \tau j^r, \quad %t_j = \tau j^r,\quad
\tau_j:=t_j-t_{j-1}%\simeq M^{-1}\,t_j^{1-1/r}
\lesssim\tau^{1/r}t_j^{1-1/r}%{\color{red}\simeq t_j/j}
\qquad %t_{j+1}\lesssim  t_{j}
\forall\,j=1,\ldots,M.
\eeq%
\end{lemma}\vspace{-0.3cm}

%Indeed,
\begin{proof}
It suffices to construct a submesh $\{ t'_k\}\subset\{t_j\}$ that satisfies \eqref{t_grid_gen}.
Let $ t'_1:=t_1$
and
$t'_k:=\min\{t_j : t_j\ge (Ck/M)^r\}$ (with an obvious modification near $t=T$). Here the constant $C$ is chosen sufficiently large
to ensure that $\{t_j\}$ is sufficiently dense within $\{t'_k\}$.
To be more precise, whenever $t_j=t'_k$ one has, in view of \eqref{t_grid_gen_gen},
$\tau_j/t_j\lesssim (\tau/t'_k)^{1/r}\lesssim C^{-1}k^{-1}$. Hence $ t'_k$ is sufficiently close to $(Ck/M)^r$, which ensures that
$\{t'_k\}$ satisfies all conditions in \eqref{t_grid_gen}.
\end{proof}

\section{Error analysis for L1-type discretizations}\label{sec_L1_error}

\subsection{Error estimation for a simplest example\! (without spatial derivatives)}\label{ssec_L1_prdgm}
It is convenient to illustrate our approach to the estimation of the temporal-discretization error using a very simple example.
Consider a fractional-derivative problem without spatial derivatives together with its discretization:
\begin{subequations}\label{simplest}
\begin{align}
D_t^\alpha u(t)&=f(t)&&\hspace{-1.6cm}\mbox{for}\;\;t\in(0,T],&&\hspace{-0.8cm} u(0)=u_0,
\\ \delta_t^\alpha U^m&=f(t_m)&&\hspace{-1.6cm}\mbox{for}\;\;m=1,\ldots,M,&&\hspace{-0.8cm} U^0=u_0.
\end{align}
\end{subequations}
Throughout this subsection, with slight abuse of notation, $\pt_t $ will be used for $\frac{d}{dt}$, while
$\delta_tu(t_j):=\tau_j^{-1}[u(t_j)-u(t_{j-1})]$.

The main result here is the following error estimate, to the proof of which we shall devote the remainder of the subsection.
\smallskip

\begin{theorem}\label{lem_simplest_star}
Let the temporal mesh
either
satisfy \eqref{t_grid_gen} with $r\ge 1$,
or
include a submesh of type \eqref{t_grid_gen}  with
the same first mesh interval.
Suppose that $u$ and $\{U^m\}$ satisfy \eqref{simplest}, and
$|\partial_t^l u |\lesssim 1+t^{\alpha-l}$ for $l = 1,2$ and $t\in(0,T]$.
Then $\forall\,m\ge1$
\beq\label{E_cal_m}
|u(t_m)-U^m|\lesssim {\mathcal E}^m:=
\left\{\begin{array}{ll}
M^{-r}\,t_m^{\alpha-1}&\mbox{if~}1\le r<2-\alpha,\\[0.3cm]
M^{\alpha-2}\,t_m^{\alpha-1}[1+\ln(t_m/t_1)]&\mbox{if~}r=2-\alpha,\\[0.2cm]
M^{\alpha-2}\,t_m^{\alpha-(2-\alpha)/r}%{\min\{0,\alpha-(2-\alpha)/r\}}
&\mbox{if~}r>2-\alpha.
\end{array}\right.
\eeq
\end{theorem}%\vspace{-0.3cm}
%\smallskip

\begin{remark}[Convergence in positive time]\label{rem_positive_time}
Consider $t_m\gtrsim 1$. Then
${\mathcal E}^m\simeq M^{-r}$ for $r<2-\alpha$ and ${\mathcal E}^{m}\simeq M^{\alpha-2}$ for $r>2-\alpha$,
i.e. in the latter case the optimal convergence rate is attained.
For $r=2-\alpha$ one gets an almost optimal convergence rate as now ${\mathcal E}^{m}\simeq M^{\alpha-2}\ln M$.

By contrast, \cite[Theorem~3.1]{sinum18_liao_et_al}
(obtained by means of a discrete Gr\"{o}nwall inequality)
 gives a somewhat similar, but less sharp error bound for graded meshes, as (in our notation)
it involves the term $O(\tau^\alpha)=O(M^{-\alpha r})$, so, e.g.,
for $r=2-\alpha$ the error bound \cite[(3.17)]{sinum18_liao_et_al} %requires (in our notation) $r=(2-\alpha)/\alpha$ to attain the optimal convergence rate in positive time.
gives a considerably less sharp convergence rate of only $\alpha(2-\alpha)$.
For $r=1$, we have ${\mathcal E}^m\simeq M^{-1}$, so our error bound is consistent with \cite{gracia_etal_cmame,laz_L1,NK_MC_L1} and is again sharper than \cite[(3.17)]{sinum18_liao_et_al}.
\end{remark}
\smallskip

\begin{remark}[Global convergence]\label{rem_global_time}
Note that $%\displaystyle
\max_{m\ge1}{\mathcal E}^m\simeq {\mathcal E}^1\simeq \tau_1^\alpha\simeq M^{-\alpha r}$ for $\alpha\le (2-\alpha)/r$,
while
$\max_{m\ge1}{\mathcal E}^m\simeq {\mathcal E}^M\simeq M^{\alpha-2}$ otherwise.
Consequently, Theorem~\ref{lem_simplest_star} yields the global error bound
$|u(t_m)-U^m|\lesssim M^{-\min \{\alpha r,2-\alpha\}} $.
This implies that
the optimal grading parameter for global accuracy is $r=(2-\alpha)/\alpha$.
Note that similar global error bounds were obtained in \cite{sinum18_liao_et_al,NK_MC_L1,stynes_etal_sinum17}.
\end{remark}
\smallskip

We first prove an auxiliary result.\smallskip

\begin{lemma}[Truncation error]\label{lem_trunc}
For a sufficiently smooth $u$, let $r^m:=\delta_{t}^{\alpha} u(t_m)-D_t^\alpha u(t_m)$ $\forall\,m\ge 1$, and
\begin{subequations}\label{psi_def}
\begin{align}\label{psi_1_def}
\psi^1&:=\sup_{s\in(0,t_1)}\!\!\bigl(s^{1-\alpha}|\delta_tu(t_1)-\pt_s u(s)|\bigr),
\\\label{psi_j_def}
\psi^j&:= t_j^{2-\alpha}\! \sup_{s\in(t_{j-1},t_j)}\!\!\!|\pt_s^2u(s)|\qquad\forall\,j\ge 2.
\end{align}
\end{subequations}
Then, under conditions \eqref{t_grid_gen} on the temporal mesh,
\beq\label{tr_er_bound}
|r^m|\lesssim (\tau/t_m)^{\min\{\alpha+1,\,(2-\alpha)/r\}} \max_{j=1,\ldots,m}\bigl\{ \psi^j\bigr\}.
\eeq
\end{lemma}\vspace{-0.3cm}

\begin{proof}
To a large degree we shall follow the proofs of \cite[Lemmas~2.3 and~2.3${}^*$]{NK_MC_L1}, so some details will be skipped.
First, recalling the definitions \eqref{CaputoEquiv} and \eqref{delta_def} of $D^\alpha_t$ and $\delta_t^\alpha$
and using the auxiliary function $\chi:=u-u^I$, we arrive at
$$
\Gamma(1-\alpha)\,r^m\!=\!\sum_{j=1}^m\!\int_{t_{j-1}}^{t_j}\!\!\!\!(t_m-s)^{-\alpha}\!\underbrace{[\delta_t u(t_j)-\pt_s u(s)]}_{{}=-\chi'(s)}ds
=\alpha%\sum_{j=1}^m\!\int_{t_{j-1}}^{t_j}
\int_0^{t_m}
\!\!\!\!(t_m-s)^{-\alpha-1}\chi(s)\,ds.
$$
On $(0,t_1)$ note that $\chi(s)=-\int_s^{t_1}\chi'(\zeta)d\zeta$, so
$|\chi(s)|\lesssim s^{\alpha-1}(t_1-s)\psi^1$  (see \cite[(2.7b)]{NK_MC_L1} for details). Otherwise,
$|\chi|\lesssim \tau_j^2t_j^{\alpha-2}\psi^j$ on $(t_{j-1},t_j)$ for $1< j<m$ and
$|\chi|\lesssim \tau_m(t_m-s)t_m^{\alpha-2}\psi^m$ on on $(t_{m-1},t_m)$.
 Consequently, a calculation shows that
\beq\label{r_m_simplest}
|r^m|\lesssim\mathring{\mathcal J}^m\,(\tau_1/t_m)^{\alpha+1}\,\psi^1+
{\mathcal J}^m\max_{j=2,\ldots,m}\bigl\{\nu_{m,j}(\tau_j/t_j)^{2-\alpha} (t_j/t_m)^{\alpha+1}%^{1-\alpha/r}
\,\psi^j\bigr\}.
\eeq
Note that in various places here we also used $  t_{j-1}\simeq t_{j}\simeq s$ for $s\in(t_{j-1},t_{j})$, $j> 1$.
The notation in~\eqref{r_m_simplest} is as follows:
%(with the use, when deriving ${\mathcal J}^m$, of $1\le t_j^{\alpha/r}s^{-\alpha/r}$ for $s\in(t_{j-1},t_j)$)
\begin{align*}
\mathring{\mathcal J}^m&:= (t_m/\tau_1)^{\alpha+1}\int_0^{t_1}\!\!s^{\alpha-1}
(t_1-s)\,%\min\{s,\,t_1-s\}
(t_m-s)^{-\alpha-1}
%\,
ds
\lesssim 1
,
%\qquad\nu^j:=(\tau_j/\tau_m)^\alpha(t_j/t_m)^{-\alpha(1-1/r)} \lesssim 1
\\
{\mathcal J}^m&:=\tau_m^\alpha\, t_m^{\alpha/r+1}\int_{t_1}^{t_m}\!
s^{-\alpha/r-1}
\,(t_m-s)^{-\alpha-1}\,\min\{1,(t_m-s)/\tau_m\}\,ds \lesssim 1
,
\\[0.2cm]
\nu_{m,j}&:=(\tau_j/\tau_m)^\alpha\,(t_j/t_m)^{-\alpha(1-1/r)} \simeq 1. %\quad\mbox{for}\;\;j=2,\ldots, m
\end{align*}
Here the bound on $\nu_{m,j}$ follows from $\tau_j/\tau_m\simeq (t_j/t_m)^{1-1/r}$ (in view of \eqref{t_grid_gen}).

For the estimation of quantities of type $\mathring{\mathcal J}^m$ and ${\mathcal J}^m$, we refer the reader to \cite{NK_MC_L1}.
In particular,
for $\mathring{\mathcal J}^m$, we first use the observation that $(t_1-s)/(t_m-s)\le t_1/t_m$ for $s\in(0,t_1)$.
Then for  $\mathring{\mathcal J}^m$ and ${\mathcal J}^m$,
  it is helpful to respectively use the substitutions $\hat s=s/t_1$ and $\hat s=s/t_m$,
while for ${\mathcal J}^m$ we also employ $(t_1/t_m)^{-\alpha/r}\sim (\tau_m/t_m)^{-\alpha}$ (also in view of \eqref{t_grid_gen}).

Combining the above observations  with \eqref{r_m_simplest}  yields
$$
|r^m|\lesssim \max_{j=1,\ldots,m}\bigl\{\underbrace{(\tau_j/t_j)^{2-\alpha}}_{\simeq(\tau/t_j)^{(2-\alpha)/r}} (t_j/t_m)^{\alpha+1}
\,\psi^j\bigr\},
$$
where we also used $\tau_j/t_j\simeq(\tau/t_j)^{1/r}$ (in view of \eqref{t_grid_gen}).
The desired bound \eqref{tr_er_bound} follows as $\tau\le t_j\le t_m$.
\end{proof}
\smallskip

\begin{corollary}[More general meshes]\label{rem_trunc_gen_mesh}
Lemma~\ref{lem_trunc} remains valid if
the temporal mesh is obtained by adding new  nodes to any mesh of type
 \eqref{t_grid_gen} under the condition that the first mesh interval remains unchanged.
\end{corollary}
\smallskip
\begin{proof}
Suppose the temporal mesh $\{t_k'\}$ is obtained by refining
the mesh $\{t_j\}_{j=0}^M$ of type \eqref{t_grid_gen}.
We again employ $\chi=u-u^I$, only now $u^I$ denotes the piecewise-linear interpolant on the finer mesh $\{t_k'\}$.
As  $t_1'=t_1$,
the estimation of integrals over $(0,t_1)$ remains unchanged.
If $t_k'=t_m$ for some $m>1$, then we  proceed exactly as in the proof of Lemma~\ref{lem_trunc},
as the same bound $|\chi|\lesssim \tau_j^2\min\{1,(t_m-s)/\tau_m\}\,t_j^{\alpha-2}\psi^j$ holds true on $(t_{j-1},t_j)$ $\forall\,j>1$ (even though $u^I$ is now the interpolant on a finer mesh).
We also use $\max_{j\le m}\psi^j\simeq \Psi'_k:=\max_{l\le k}\{\psi'^l\}$.
If $t_k'\in(t_{m-1},t_m)$ for some $m>1$, then one has a similar bound $|\chi|\lesssim \tau_j^2\min\{1,(t_k'-s)/\tau_m\}\,t_j^{\alpha-2}\,\Psi'_k$ on $(t_{j-1},\min\{t_j,t_k'\})$ $\forall\,j>1$.
So for  the truncation error
at $t_k'$ we get a version of \eqref{r_m_simplest}, in which (including its ingredients)
 $t_m$ is replaced by $t_k'\simeq t_m$ and $\psi^j$ is replaced by $\Psi'_k$.
 This again leads to the desired version of \eqref{tr_er_bound} for the mesh $\{t_k'\}$.
\end{proof}
\smallskip

{\it Proof of Theorem~\ref{lem_simplest_star}.}
%\begin{proof}
Consider the error $e^m:=u(t_m)-U^m$, for which \eqref{simplest} implies
$e^0=0$ and $\delta_t^\alpha e^m=r^m$  $\forall\,m\ge 1$, where
the truncation error $r^m$ is from Lemma~\ref{lem_trunc} and hence satisfies \eqref{tr_er_bound}.
Furthermore, under the conditions \eqref{t_grid_gen} on the temporal mesh (or its submesh),
one has $\psi^1\lesssim 1$ (in view of $|\delta_tu(t_1)|\le\tau_1^{-1}\int_0^{t_1}|\pt_s u|\,ds\lesssim \tau_1^{\alpha-1}$)
and $\psi^j\lesssim 1$ for $j\ge 2$ (in view of $s\simeq t_j$ for $s\in(t_{j-1},t_j)$ for this case).
Consequently, in view of Lemma~\ref{lem_trunc} and Corollary~\ref{rem_trunc_gen_mesh}, we arrive at
$$
|r^m|\lesssim (\tau/t_m)^{\gamma+1}\quad\forall\,m\ge 1,
\qquad\mbox{where}\;\;
\gamma+1:=\min\{\alpha+1,(2-\alpha)/r\}.
$$

Next consider three cases.
\smallskip

Case $1\le r<2-\alpha$. Then both $(2-\alpha)/r> 1$ and $\alpha+1>1$,
so $\gamma>0$.
An application of Theorem~\ref{theo_barrier}(i) for this case yields
$|e^m|\lesssim \tau\, t_m^{\alpha-1}$, where $\tau\sim M^{-r}$.
\smallskip

Case $r=2-\alpha$. Then  $(2-\alpha)/r= 1$, while $\alpha+1>1$,
so $\gamma=0$.
An application of Theorem~\ref{theo_barrier}(i) yields
$|e^m|\lesssim \tau\, t_m^{\alpha-1}[1+\ln(t_m/t_1)]$,
where $\tau\sim M^{-r}=M^{\alpha-2}$.
\smallskip

Case $r>2-\alpha$.
Then
$(2-\alpha)/r< 1$, while $\alpha+1>1$,
so $\gamma+1=(2-\alpha)/r<1$.
An application of Theorem~\ref{theo_barrier}(where part~(i) of this theorem is used if $r\le (2-\alpha)/\alpha$
and part~(ii) is used otherwise)
 yields
$|e^m|\lesssim \tau\, t_m^{\alpha-1}(\tau/t_m)^{(2-\alpha)/r-1}\sim
\tau^{(2-\alpha)/r}t_m^{\alpha-(2-\alpha)/r}$,
where
$\tau^{(2-\alpha)/r}\sim M^{\alpha-2}$.
%
%Case $2-\alpha<r\le (2-\alpha)/\alpha$.
%Then
%$(2-\alpha)/r< 1$, while $\alpha+1>1$,
%so $\gamma+1=(2-\alpha)/r\in[\alpha,1)$.
%An application of Theorem~\ref{theo_barrier}(i) for this case yields
%$|u(t_m)-U^m|\lesssim \tau\, t_m^{\alpha-1}(\tau/t_m)^{(2-\alpha)/r-1}\sim
%\tau^{(2-\alpha)/r}t_m^{\alpha-(2-\alpha)/r}$,
%where
%$\tau^{(2-\alpha)/r}\sim M^{2-\alpha}$.
%\smallskip
%
%Case $r> (2-\alpha)/\alpha$. Then $(2-\alpha)/r< \alpha<\alpha+1$,
%so $\gamma+1=(2-\alpha)/r< \alpha$.
%Hence $|r^m|\lesssim \tau^{\gamma+1}t_m^{-(\gamma+1)}\lesssim \tau^{\gamma+1}t_m^{-\alpha}$, so
%application of Theorem~\ref{theo_barrier}(ii) yields
%$|u(t_m)-U^m|\lesssim \tau^{\gamma+1}=\tau^{(2-\alpha)/r}\sim M^{2-\alpha}$.
%
%\end{proof}
\endproof

\subsection{Error analysis for the L1 semidiscretization in time}\label{ssec_L1_semi}

Consider the semidiscretization of our problem~\eqref{problem} in time using the L1 method:
%$\delta_t^\alpha$ from \eqref{delta_t_def}:
\beq\label{semediscr_method}
\delta_t^\alpha U^m +\LL U^m= f(\cdot,t_m)\;\;\mbox{in}\;\Omega,\quad U^m=0\;\;\mbox{on}\;\pt\Omega\quad\forall\,m=1,\ldots,M;\quad U^0=u_0.
\eeq

\begin{theorem}\label{theo_semidiscr}
Let the temporal mesh
%satisfy \eqref{t_grid_gen} with $r\ge 1$.
%
either
satisfy \eqref{t_grid_gen} with $r\ge 1$,
or
include a submesh of type \eqref{t_grid_gen}  with
the same first mesh interval.
 Given $p\in\{2,\infty\}$,
 suppose that $u$ is from \eqref{problem},\eqref{LL_def}, with
 \mbox{$c-p^{-1}\!\sum_{k=1}^d\!\pt_{x_k}\!b_k\ge 0$},
and %it satisfies
$\|\partial_t^l u (\cdot, t)\|_{L_p(\Omega)}\lesssim 1+t^{\alpha-l}$ for $l = 1,2$ and $t\in(0,T]$.
Then for  $\{U^m\}$ from \eqref{semediscr_method}, one has
\beq\label{L1_semi_error}
\|u(\cdot,t_m)-U^m\|_{L_p(\Omega)}\lesssim {\mathcal E}^m
\qquad\forall\,m=1,\ldots,M,
\eeq
where ${\mathcal E}^m$ is from \eqref{E_cal_m}.
\end{theorem}
\smallskip

\begin{proof}
For the error $e^m:= u(\cdot,t_m)-U^m$, using \eqref{problem} and \eqref{semediscr_method}, and imitating the proof of \cite[Theorem~3.1]{NK_MC_L1}, one gets a version of
\cite[(3.4)]{NK_MC_L1}:
\beq\label{semi_error}
\delta_t^\alpha \|e^m\|_{L_p(\Omega)}\le \|r^m\|_{L_p(\Omega)}%+\vartheta\|e^m\|_{L_p(\Omega)}
\qquad\forall\,m=1,\ldots, M.
\eeq
Here
the truncation error $r^m:=\delta_t^\alpha u(\cdot,t_m)-D_t^\alpha u(\cdot,t_m)$ is estimated in Lemma~\ref{lem_trunc} and hence satisfies \eqref{tr_er_bound}.
The desired error bound is obtained by closely imitating the proof of Theorem~\ref{lem_simplest_star}. Importantly, parts (i) and (ii) of Theorem~\ref{theo_barrier}
remain applicable to \eqref{semi_error} in view of Theorem~\ref{theo_barrier}(iii).
\end{proof}

\subsection{Error analysis for full L1-type discretizations}\label{ssec_L1_full}

Similarly to \S\ref{ssec_L1_semi}, one can easily combine the analysis of \S\ref{ssec_L1_prdgm} with \cite[\S\S4-5]{NK_MC_L1}
to obtain error bounds of type \eqref{E_cal_m}
for full discretizations of problem \eqref{problem} with $\LL=\LL(t)$, whether finite differences or finite elements are employed as spatial discreziations.
We shall give a flavour of such results.
%\smallskip

\subsubsection{Finite difference discretizations}\label{ssec_FD}
Consider our problem \eqref{problem}--\eqref{LL_def} in the spatial domain $\Omega=(0,1)^d\subset\R^d$.
Suppose that the standard finite difference operator $\LL_h $ from \cite[\S4]{NK_MC_L1} is employed
as a spatial discretization on a uniform tensor-product mesh $\Omega_h$ of size $h$.
We shall assume that $h$ is sufficiently small so that $\LL_h $ satisfies the discrete maximum principle.
Then, under the conditions of Theorem~\ref{theo_semidiscr} with $p=\infty$, and additionally assuming that
$\|\pt_{x_k}^l u(\cdot,t)\|_{L_\infty(\Omega)}\lesssim 1$ for $l=3,4$, $k=1,\ldots,d$ and $t\in(0,T]$, one easily gets
the following version of \cite[Theorem~4.1]{NK_MC_L1}:
\beq\label{FD_error}
\|u(\cdot,t_m)-U^m\|_{\infty\,;\Omega_h}\lesssim {\mathcal E}^m+t_m^{\alpha}\, h^2
\qquad\forall\,m=1,\ldots,M,
\eeq
where
$\|\cdot\|_{\infty\,;\Omega_h}:=\max_{\Omega_h}|\cdot|$
denotes the spatial nodal maximum norm, while ${\mathcal E}^m$ is from \eqref{E_cal_m}.

\subsubsection{Finite element discretizations}\label{ssec_L1_full_FEM}
Discretize \eqref{problem}--\eqref{LL_def}, posed in a general bounded  Lipschitz domain  $\Omega\subset\R^d$,
 by applying
 a standard Galerkin finite element spatial approximation to the temporal semidiscretization~\eqref{semediscr_method}.
  A Lagrange finite element space $S_h \subset H_0^1(\Omega)\cap C(\bar\Omega)$ of fixed degree $\ell\ge 1$, % \Bbbk
 relative to a
 quasiuniform simplicial triangulation
 %{\color{red}shape regular} mesh
 %$\mathcal T$
 of $\Omega$, is employed, as in \cite[\S5]{NK_MC_L1}.
 Then,  under the conditions of Theorem~\ref{theo_semidiscr} with $p=2$,
% and additionally assuming that $$
%$\|\pt_{x_k}^l u(\cdot,t)\|_{L_\infty(\Omega)}\lesssim 1$ for $l=3,4$, $k=1,\ldots,d$ and $t\in(0,T]$,
one easily gets
the following version of \cite[Theorem~5.1]{NK_MC_L1}:
\begin{align}\notag
\|u(\cdot,t_m)-u_h^m\|_{L_2(\Omega)}\lesssim{}&\|u_0-u_h^0\|_{L_2(\Omega)}+{\mathcal E}^m
+\max_{t\in\{0,t_m\}}\|\rho(\cdot, t)\|_{L_2(\Omega)}
\\&{}+\int_0^{t_m}\|\pt_t \rho(\cdot, t)\|_{L_2(\Omega)}\,dt
\qquad\forall\,m=1,\ldots,M.
\label{L1_fem_error}
\end{align}
Here $u_h^m\in S_h$ is the finite element solution at time $t_m$, ${\mathcal E}^m$ is from \eqref{E_cal_m}, and
$\rho(\cdot, t):=\RR_h u(t)-u(\cdot, t)$ is the error of
the standard Ritz projection $\RR_h u(t)\in S_h$ of $u(\cdot,t)$.
Under additional realistic assumptions on $u$, the final two terms in the above error estimate can be bounded by $O( h^{\ell+1})$,
where $h$ is the triangulation diameter \cite[\S5]{NK_MC_L1}.
One can also estimate the error in the norm $L_\infty(\Omega)$ imitating \cite[\S5.2]{NK_MC_L1}.

\section{Generalization for the Alikhanov discrete fractional-derivative operator}\label{sec_A_stab}
In this section we shall show that the above error analysis is not restricted to L1 discretizations, but may be extended, without major modifications, to other discretizations.
Here we shall focus on a higher-order discrete fractional-derivative operator proposed by Alikhanov \cite{alikh},
while  a similar analysis is generalized for another higher-order scheme  in \cite{NK_L2}.

A stability property of type \eqref{main_stab} will be established in \S\ref{sssec_Alikh_B_lemma}. Next, in \S\ref{ssec_A_prdgm}, the truncation error will be estimated and the error for the simplest problem without spatial derivatives
will be bounded by a quantity similar to ${\mathcal E}^m$  in \eqref{E_cal_m}.
A stability property of type \eqref{main_stab} for Alikhanov-type semi-discretizations will be obtained in~\S\ref{ssec_A_semi}, which will allow to extend our error analysis to this case.
Finally, error bounds for full discretizations will be briefly discussed in \S\ref{ssec_A_full}.

\subsection{Alikhanov discrete fractional-derivative operator. Discrete maximum principle}
The discrete fractional-derivative operator proposed by Alikhanov is associated with the point
\begin{subequations}\label{delta_t_def_ab}
\beq\label{t_star}\textstyle
t_m^*:=t_{m-\alpha/2}=t_m-\frac12\alpha\tau_m.
\eeq
In the definition of this operator, as well as in its analysis, we shall employ three standard
{\it Lagrange interpolation operators} with the following interpolation points:
$$
\Pi_{1,j}\;:\;\{t_{j-1},t_j\},\qquad
\Pi_{2,j}\;:\;\{t_{j-1},t_j,t_{j+1}\},\qquad
\Pi^*_{2,j}\;:\;\{t_{j-1},t_j^*,t_j\}
$$
Now, applying $\Pi_{2,j}$ to the computed solution values $\{U^j\}$ on $(t_{j-1},t_j)$ for $j<m$
and $\Pi_{1,m}$ on on $(t_{m-1},t_m^*)$,
we define an alternative discretization for the fractional operator $D_t^\alpha$ $\forall\,m= 1,\ldots, M$:
\beq\label{delta_t_def}
\delta_{t}^{\alpha,*} U^m :=D^\alpha_t (\Pi^mU)(t^*_m),\quad
\quad%\mbox{where}\quad
\Pi^m:=\left\{\begin{array}{cl}
%\Pi_{1,1}&\mbox{on~}(0,t_1)&\mbox{for~}m=1,\\
\Pi_{2,j}&\mbox{on~}(t_{j-1},t_j)\;\;\; \forall j<m,\\
\Pi_{1,m}&\mbox{on~}(t_{m-1},t^*_m).\\
\end{array}\right.
\eeq
\end{subequations}
Note that the interpolation operator $\Pi^*_{2,j}$ is not used in the definition of $\delta_{t}^{\alpha,*}$, but will be useful
in the estimation of the truncation error.
In particular, for the final interval $(t_{m-1},t_m^*)$ it will occasionally be convenient
to employ the representation $\Pi_{1,m}=\Pi^*_{2,m}+(\Pi_{1,m}-\Pi^*_{2,m})$, as
the choice \eqref{t_star} ensures
for any sufficiently smooth function $v$ that
\beq\label{nu_m}
\int_{t_{m-1}}^{ t_m^*}\!\!(\Pi_{1,m}v-\Pi^*_{2,m}v)'(s)\,(t^*_m-s)^{-\alpha}ds=0.
\eeq
Indeed, here
$\Pi_{1,m}v-\Pi^*_{2,m}v=C(s-t_{m-1})(t_m-s)$, with some constant $C$,
so one has
$(\Pi_{1,m}v-\Pi^*_{2,m}v)'=2C(s-t_{m-1/2})$ and, consequently,  \eqref{t_star} yields \eqref{nu_m}.

\smallskip

\begin{remark}[Discrete maximum principle]\label{rem_JSC}
Sufficient conditions for the operator $\delta_{t}^{\alpha,*}$
to be associated with an M-matrix, and, hence, satisfy the discrete maximum principle, are given by \cite[Lemma~4]{ChenMS_JSC} (see also \cite[Remark~3]{ChenMS_JSC}) and \cite[assumption M1]{arxiv18_liao_et_al}.
In particular, throughout this section we shall
assume that either ${0.4656\le{}}\rho_{j}\le\rho_{j-1}$ $\forall\,j\ge 2$ \cite{ChenMS_JSC}
or $\rho_j\ge 4/7$
$\forall\,j\ge 2$ \cite{arxiv18_liao_et_al}, where $\rho_j:=\tau_{j+1}/\tau_j$.
It is sufficient for the discrete maximum principle, and it is satisfied, for example, by the standard graded mesh
$\{t_j=T(j/M)^r\}_{j=0}^M$ with any $r\ge 1$.
\end{remark}

\subsection{Stability theorem for
the Alikhanov scheme}\label{sssec_Alikh_B_lemma}

To generalize the above error analysis to the Alikhanov scheme, we need to extend the stability result given by Theorem~\ref{theo_barrier} to the operator $\delta_{t}^{\alpha,*}$.
\smallskip

\begin{theorem}[Stability]\label{theo_barrierA}
Let the temporal mesh $\{t_j\}_{j=0}^M$  satisfy the condition from Remark~\ref{rem_JSC}.\\[0.1cm]
(i) Let the temporal mesh additionally satisfy
 \eqref{t_grid_gen}
with $1\le r\le (3-\alpha)/\alpha$.
Given  $\gamma\in \R$ and $\{V^j\}_{j=0}^M$,
% with $V^0=0$ one has
the stability property \eqref{main_stab}  holds true with $\delta_{t}^{\alpha}$ replaced by $\delta_{t}^{\alpha,*}$.
%for $j=1,\ldots, M$.\\
\\[0.1cm]
(ii)
If $\gamma\le \alpha-1$, then, without further restrictions on the mesh, \eqref{main_stab} holds true   with $\delta_{t}^{\alpha}$ replaced by $\delta_{t}^{\alpha,*}$.
\\[0.1cm]
(iii) The above results remain valid if  $|\delta_{t}^{\alpha,*} V^j|\lesssim (\tau/ t_j)^{\gamma+1}$ in \eqref{main_stab} is replaced by $\delta_{t}^{\alpha,*} |V^j|\lesssim (\tau/ t_j)^{\gamma+1}$.%
\end{theorem}
\smallskip

\begin{proof}
(i)
This part is obtained similarly to the proof of Theorem~\ref{theo_barrier}(i), only with a few changes in obtaining
 a version of Lemma~\ref{lem_stability1} for $\delta_{t}^{\alpha,*}$; see Lemma~\ref{lem_stability1_A} below.
\smallskip

(ii)
This part is obtained exactly as in the proof of Theorem~\ref{theo_barrier}, only instead of
 \cite[Lemma 2.1(i)]{NK_MC_L1} we now employ a similar \cite[Lemma~5]{ChenMS_JSC} for $\delta_{t}^{\alpha,*}$.
\smallskip

(iii) This part is obtained exactly as in the proof of Theorem~\ref{theo_barrier}.
\end{proof}
\smallskip

\begin{lemma}[Lemma~\ref{lem_stability1} for
Alikhanov scheme]\label{lem_stability1_A}
Under the conditions of Theorem~\ref{theo_barrierA}(i) on the temporal mesh, the
 discrete barrier function $\{B^j\}_{j=0}^M$  from
 \eqref{B_L1_def}
 %Lemma~\ref{lem_stability1}
 satisfies
\eqref{B_def} with $\delta_{t}^{\alpha}$ replaced by $\delta_{t}^{\alpha,*}$.
\end{lemma}
\smallskip

\begin{proof}
As $t_j^*\sim t_j$ (in view of \eqref{t_star}), it suffices to prove that $\delta_{t}^{\alpha,*} B^j\ge \tau^\alpha (t_j^*)^{-\alpha-1}$ $\forall j\ge1$.
For the latter, we closely imitate the proof of Lemma~\ref{lem_stability1}.
In particular, for $j\le p$ one gets $\delta_{t}^{\alpha,*} B^j=D_t^\alpha B(t_j^*)$.
When estimating $\delta_{t}^{\alpha,*} B^m$ for  $m>p$,
 a few modifications are required
that we now describe.

For $D^\alpha_t B(t_m^*)$ we have \eqref{app_B_eq}, while,
in view of \eqref{nu_m},
$\delta_{t}^{\alpha,*} B^m=D^\alpha_t (I_2B)(t_m^*)$,
where
$I_2B:=\Pi_{2,j}B$ on $(t_{j-1},t_j)$ for $j<m$ and $I_2B:=\Pi^*_{2,j}B$ on $(t_{m-1},t_m^*)$ (with interpolation points $\{t_{m-1},t_m^*,t_m\}$),
i.e. $I_2B$ is a piecewise quadratic interpolant.
Now
$\Gamma(1-\alpha)[\delta_{t}^{\alpha,*} B^m-D^\alpha_t B(t_m^*)]=\sum_{j=p}^m\mu^j$, where
(compare with \eqref{mu_def})
\begin{align*}
\mu^j&:= %\int_{t_{j-1}}^{t_j}\!\![\delta_t B^j-B'(s)]\underbrace{(t_m-s)^{-\alpha}}_{{}=g(s)}ds=
%\int_{t_{j-1}}^{t_j}\!\!(B^I-B)'(s)\,(t_m-s)^{-\alpha}ds
\alpha\int_{t_{j-1}}^{\min\{t_j,t_m^*\}}\!\! (B-I_2B)(s)\,(t^*_m-s)^{-\alpha-1} ds.
%\\
%\nu^m&:=\int_{t_{m-1}}^{ t_m^*}\!\!(\Pi_{1,m}B-\Pi^*_{2,m}B)'(s)\,(t^*_m-s)^{-\alpha}ds.
\end{align*}
%Importantly, here we employ
%
%For $\nu^m$, note that
%$\Pi_{1,m}B-\Pi^*_{2,m}B=C(s-t_{m-1})(t_m-s)$, with some constant $C$,
%so
%$(\Pi_{1,m}B-\Pi^*_{2,m}B)'=2C(s-t_{m-1/2})$ and, consequently,  \eqref{t_star} yields $\nu^m=0$.

The estimation of $\mu^j$ for $j>p$ is similar to the case of the L1 scheme,
only now we use a sharper bound
$| B-I_2 B |\lesssim \tau_j^3 \min\{1,(t_m^*-s)/(t_m^*-t_{m-1})\}|B'''(t_{j-1})|$, where $|B'''(t_{j-1})|\lesssim s^{-\beta-3}$.
So now we get the following version of
\eqref{mu_sum}, in which the factors that differ from the proof of Lemma~\ref{lem_stability1} are framed:
$$
\sum_{p+1}^m|\mu^j|\lesssim
 \tau^{\fbox{\scriptsize$3$}/r}\!\int_{t_{p}}^{t^*_m}\!\!\!s^{-\beta-\fbox{\scriptsize$3$}/r}\,(t_m-s)^{-\alpha-1}\,\min\bigl\{1,(t^*_m-s)/(t_m^*-t_{m-1})\bigr\}\,ds.
$$
This leads to the following version of \eqref{mu_sum1}:
\begin{align}\notag
\sum_{p+1}^m|\mu^j|&\lesssim  \tau^{\fbox{\scriptsize$3$}/r}(t^*_m)^{-\fbox{\scriptsize$3$}/r-1}\Bigl[ (t_p/t_m^*)^{\alpha-\fbox{\scriptsize$3$}/r}  +(\tau/t^*_m)^{-\alpha/r}\Bigr]\\[0.2cm]\notag
&\lesssim  (\tau/t_p)^{\fbox{\scriptsize$3$}/r}\, t_p^\alpha (t^*_m)^{-\alpha-1}
+ \underbrace{(t^*_m)^{-1}(\tau/t^*_m)^{(\fbox{\scriptsize$3$}-\alpha)/r}}_{\lesssim \tau^\alpha (t_m^*)^{-\alpha-1}}\\
\label{sum_mu_Alikh}
&\lesssim
\Bigl[(\tau/t_p)^{\fbox{\scriptsize$3$}/r}+(\tau/t_p)^\alpha\Bigr]
\,t_p^\alpha (t^*_m)^{-\alpha-1},
\end{align}
where in the second line we employed
$(\tau/t^*_m)^{(3-\alpha)/r}\lesssim (\tau/t^*_m)^{\alpha}$ (in view of $r\le (3-\alpha)/\alpha$).

It remains to get a similar bound on $|\mu^p|$ (where $p<m$).
As $B'$ abruptly changes at $t_p$, we now employ
$| B-I_2 B |=| B-\Pi_{2,j} B |\lesssim \max_{[t_{p-1},t_{p+1}]}|B-B(t_p)|\lesssim \tau_p t_p^{-\beta-1}$.
(Note that when using the latter bound, we rely on the property $\tau_j\sim\tau_{j+1}$ for the stability of the interpolating operator $\Pi_{2,j}$
in the sense that
$\max_{[t_{j-1},t_{j}]}|\Pi_{2,j} v|\lesssim \max_{[t_{j-1},t_{j+1}]}|v|$ for any continuous $v$.)
Now a calculation shows that
$$
|\mu^p|%\lesssim \int_{t_{p-1}}^{t_p}\!\! [\Pi_{2,j}B-B](\tilde t_m-s)^{-\alpha-1} ds
\lesssim \tau_p t_p^{-\beta-1}\int_{t_{p-1}}^{t_p}\!\! (t^*_m-s)^{-\alpha-1} ds
\lesssim \tau_p t_p^{-\beta-1}(t^*_{p+1}/ t^*_m)^{\alpha+1} \underbrace{\int_{t_{p-1}}^{t_p}\!\! (t^*_{p+1}-s)^{-\alpha-1} ds}_{\lesssim \tau_p^{-\alpha}}\,,
$$
where, in the second relation, we employed the observation $(t^*_{p+1}-s)/ (t^*_m-s)\le t^*_{p+1}/ t^*_m$ $\forall s\in(0,t^*_{p+1})$
(in view of $t^*_{p+1}\le t^*_m$).
Next,  $|\mu^p|\lesssim
 (\tau_p/t_p)^{\beta}\,t_{p}^{\alpha}  t_m^{-\alpha-1}$, and, in view of
 $\tau_p/t_p=(\tau/t_p)^{1/r}$ (by \eqref{t_grid_gen} ), one gets
$|\mu^p|\lesssim (\tau/t_p)^{\beta/r}\,t_{p}^{\alpha}  t_m^{-\alpha-1}$.

Finally, combining the latter bound with \eqref{sum_mu_Alikh},
we conclude that
$|\delta_{t}^{\alpha,*}\, B^m-D^\alpha_t B(t_m^*)|\lesssim \sum_{j=p}^m|\mu^j|$ will be dominated by $\frac12D^\alpha_t B(t_m^*)$ from \eqref{app_B_eq}
if $p$ is chosen sufficiently large.
\end{proof}

\subsection{Error analysis of the Alikhanov scheme for a simplest example\! (without spatial derivatives)}\label{ssec_A_prdgm}

Consider a fractional-derivative problem without spatial derivatives together with its discretization using $\delta_{t}^{\alpha,*}$ from \eqref{delta_t_def_ab}:
\begin{subequations}\label{simplestA}
\begin{align}
D_t^\alpha u(t)&=f(t)&&\hspace{-1.6cm}\mbox{for}\;\;t\in(0,T],&&\hspace{-0.8cm} u(0)=u_0,
\\ \delta_{t}^{\alpha,*} U^m&=f(t_m^*)&&\hspace{-1.6cm}\mbox{for}\;\;m=1,\ldots,M,&&\hspace{-0.8cm} U^0=u_0.
\end{align}
\end{subequations}
Then for the error we have a version of Theorem~\ref{lem_simplest_star}.
\smallskip

\begin{theorem}\label{lem_simplest_starA}
Let the temporal mesh satisfy the condition from Remark~\ref{rem_JSC} and
 \eqref{t_grid_gen}
with $r\ge 1$.
Suppose that $u$ and $\{U^m\}$ satisfy \eqref{simplest}, and
$|\partial_t^l u |\lesssim 1+t^{\alpha-l}$ for $l = 1,3$ and $t\in(0,T]$.
Then $\forall\,m\ge1$
\beq\label{E_cal_mA}
|u(t_m)-U^m|\lesssim {\mathcal E}^{m,*}:=
\left\{\begin{array}{ll}
M^{-r}\,t_m^{\alpha-1}&\mbox{if~}1\le r<3-\alpha,\\[0.3cm]
M^{\alpha-3}\,t_m^{\alpha-1}[1+\ln(t_m/t_1)]&\mbox{if~}r=3-\alpha,\\[0.2cm]
M^{\alpha-3}\,t_m^{\alpha-(3-\alpha)/r}%{\min\{0,\alpha-(2-\alpha)/r\}}
&\mbox{if~}r>3-\alpha.
\end{array}\right.
\eeq
\end{theorem}

\begin{remark}[Convergence in positive time]\label{rem_positive_time_A}
Consider $t_m\gtrsim 1$. Then
${\mathcal E}^{m,*}\simeq M^{-r}$ for $r<3-\alpha$ and ${\mathcal E}^{m,*}\simeq M^{\alpha-3}$ for $r>3-\alpha$,
i.e. in the latter case the optimal convergence rate is attained.
For $r=3-\alpha$ one gets an almost optimal convergence rate as now ${\mathcal E}^{m,*}\simeq M^{\alpha-3}\ln M$.
%
%By contrast, \cite[Theorem~3.9]{arxiv18_liao_et_al} (obtained by means of a discrete Gr\"{o}nwall inequality \cite{sinum19_liao_et_al}) gives a somewhat similar, but less sharp error bound for graded meshes, as (in our notation)
%it involves the term $O(\tau^\alpha)=O(M^{-\alpha r})$, so, e.g.,
%it  requires $r=(3-\alpha)/\alpha$ to attain the optimal convergence rate in positive time.
%For $r=1$, we have ${\mathcal E}^{m,*}\simeq M^{-1}$, so our error bound  is also sharper than those in \cite[Theorem~3.9]{arxiv18_liao_et_al}.%
\end{remark}%
\smallskip

\begin{remark}[Global convergence]\label{rem_global_time_A}
Note that $%\displaystyle
\max_{m\ge1}{\mathcal E}^{m,*}\lesssim {\mathcal E}^{1,*}\simeq \tau_1^\alpha\simeq M^{-\alpha r}$ for $\alpha\le (3-\alpha)/r$,
while
$\max_{m\ge1}{\mathcal E}^{m,*}\simeq {\mathcal E}^{M,*}\simeq M^{\alpha-3}$ otherwise.
Consequently, Theorem~\ref{lem_simplest_starA} yields the global error bound
$|u(t_m)-U^m|\lesssim M^{-\min \{\alpha r,3-\alpha\}} $.
This implies that
the optimal grading parameter for global accuracy is $r=(3-\alpha)/\alpha$.
Note that a similar global error bound was obtained in \cite{ChenMS_JSC}.
\end{remark}
\smallskip

The proof is, to a large degree, similar to the arguments in \S\ref{ssec_L1_prdgm}, with slight modifications in the truncation error estimation.
\smallskip

\begin{lemma}[Truncation error]\label{lem_truncA}
For a sufficiently smooth $u$, let $r^m:=\delta_{t}^{\alpha,*} u(t_m)-D_t^\alpha u(t^*_m)$ $\forall\,m\ge 1$, and
\begin{subequations}\label{psi_defA}
\begin{align}\label{psi_1_defA}
\psi^1&:=\sup_{s\in(0,t_2)}\!\!\bigl(s^{1-\alpha}|\pt_s u(s)|\bigr)+t_2^{-\alpha}{\rm osc}\bigl(u, [0,t_2]\bigr),
\\\label{psi_j_defA}
\psi^j&:= t_j^{3-\alpha}\!\! \sup_{s\in(t_{j-1},t_{j+1})}\!\!\!|\pt_s^3u(s)|\quad\forall\,2\le j\le M-1,\qquad \psi^M:=\psi^{M-1}.
\end{align}
\end{subequations}
Then, under conditions \eqref{t_grid_gen} on the temporal mesh, one has
\beq\label{tr_er_boundA}
|r^m|\lesssim (\tau_1/t_m)^{\min\{\alpha+1,\,  (3-\alpha)/r\}} \max_{j=1,\ldots,m}\bigl\{ \psi^{j}\bigr\}
\qquad\forall\,m\ge 1.
\eeq
\end{lemma}\vspace{-0.3cm}

\begin{proof}
We imitate the proof of Lemma~\ref{lem_trunc}, and also use the notation $I_2$ and some observations from the proof of Lemma~\ref{lem_stability1_A}.
Recall that,
in view of \eqref{nu_m},
$\delta_{t}^{\alpha,*} u(t_m)=D^\alpha_t (I_2u)(t_m^*)$
where
$I_2=\Pi_{2,j}$ on $(t_{j-1},t_j)$ for $j<m$ and $I_2:=\Pi^*_{2,j}$ on $(t_{m-1},t_m^*)$.
%So split $\delta_{t}^{\alpha,*} u(t_m)=D^\alpha_t (I_2u)(t_m^*)+\vartheta^m$, where
%$$
%\Gamma(1-\alpha)\,\vartheta^m:= \int_{t_{m-1}}^{t_m^*}\!\!(t^*_m-s)^{-\alpha}\,\pt_s (\Pi_{1,m}u-\Pi^*_{2,m}u)(s)\, ds = 0.
%$$
%Here $\vartheta^m=0$ is obtained, in view of  \eqref{t_star}, similarly to $\nu^m=0$ in the proof of the proof of Lemma~\ref{lem_trunc}.
%Hence, we conclude that $\delta_{t}^{\alpha,*} u(t_m)=D^\alpha_t (I_2u)(t_m^*)$.

Next, recalling the definition \eqref{CaputoEquiv} %and \eqref{delta_def}
of $D^\alpha_t$ %and $\delta_t^\alpha$
and using the auxiliary function $\chi:=u-I_2u$, which satisfies $\chi(t_m^*)=0$, we arrive at
$$
\Gamma(1-\alpha)\,r^m\!=\!%\sum_{j=1}^m\!\int_{t_{j-1}}^{\min\{t_j,t_m^*\}}
\int_0^{t_m^*}
\!\!\!\!(t^*_m-s)^{-\alpha}\underbrace{\pt_s[I_2u(s)- u(s)]}_{{}=-\chi'(s)}ds
=\alpha
%\sum_{j=1}^m\!\int_{t_{j-1}}^{\min\{t_j,t_m^*\}}
\!\int_0^{t_m^*}
\!\!\!\!(t^*_m-s)^{-\alpha-1}\chi(s)\,ds.
$$

Let $t_1^{**}:=\min\{t_1,t_m^*\}$
and consider the intervals $(0,t_1^{**})$ and $(t_1^{**},t_m^*)$ separately.
On $(0,t_1^{**})$ note that $\chi(t_1^{**})=0$ implies $\chi(s)=-\int_s^{t_1^{**}}\chi'(\zeta)d\zeta$,
where
$|\chi'|\le|\pt_su|+|\pt_s (I_2u)|$, while
$|\pt_s (I_2u)|\lesssim t_2^{-1}{\rm osc}(u, [0,t_2])\le s^{\alpha-1}t_2^{-\alpha}{\rm osc}(u, [0,t_2])$
(in view of $s\le \tau_1\simeq \tau_2$),
so
$|\chi(s)|\lesssim s^{\alpha-1}(t_1^{**}-s)\psi^1$.
Note also that
$|\chi|\lesssim \tau_j^3t_j^{\alpha-3}\psi^j$ on $(t_{j-1},t_j)$ for $1< j<m$ and
$|\chi|\lesssim \tau^2_m(t^*_m-s)t_m^{\alpha-3}\psi^m$  on $(t_{m-1},t^*_m)$ if $m>1$.
Consequently, a calculation shows that we get a version of \eqref{r_m_simplest}:
\beq\label{r_m_simplestA}
|r^m|\lesssim\mathring{\mathcal J}^m\,(\tau_1/t_m)^{\alpha+1}\,\psi^1+
{\mathcal J}^m\max_{j=2,\ldots,m}\bigl\{\nu_{m,j}(\tau_j/t_j)^{\fbox{\scriptsize$3$}-\alpha} (t_j/t_m)^{\alpha+1}%^{1-\alpha/r}
\,\psi^j\bigr\},
\eeq
where
for convenience, the factors that differ from the proof of Lemma~\ref{lem_trunc} are framed.
Note that in various places we also use $ t_j^*\simeq t_j\simeq t_{j+1}$ for $j\ge 1$ and $s\simeq t_j$ on $(t_{j-1},t_j)$.
The notation in \eqref{r_m_simplestA} is as follows:
%(with the use, when deriving ${\mathcal J}^m$, of $1\le t_j^{\alpha/r}s^{-\alpha/r}$ for $s\in(t_{j-1},t_j)$)
\begin{align*}
\mathring{\mathcal J}^m&:= (t_m/\tau_1)^{\alpha+1}\int_0^{t^{**}_1}\!\!s^{\alpha-1}
(t^{**}_1-s)\,%\min\{s,\,t_1-s\}
(t^*_m-s)^{-\alpha-1}
%\,
ds
\lesssim 1
,
%\qquad\nu^j:=(\tau_j/\tau_m)^\alpha(t_j/t_m)^{-\alpha(1-1/r)} \lesssim 1
\\
{\mathcal J}^m&:=\tau_m^\alpha\, t_m^{\alpha/r+1}\int_{t^{**}_1}^{t^*_m}\!
s^{-\alpha/r-1}
\,(t^*_m-s)^{-\alpha-1}\,\min\{1,(t^*_m-s)/\tau_m\}\,ds \lesssim 1
,
\\[0.2cm]
\nu_{m,j}&:=(\tau_j/\tau_m)^\alpha\,(t_j/t_m)^{-\alpha(1-1/r)} \simeq 1. %\quad\mbox{for}\;\;j=2,\ldots, m
\end{align*}
Here the bound on $\nu_{m,j}$ follows from $\tau_j/\tau_m\simeq (t_j/t_m)^{1-1/r}$ (in view of \eqref{t_grid_gen}).
For the estimation of quantities of type $\mathring{\mathcal J}^m$ and ${\mathcal J}^m$, we refer the reader to \cite{NK_MC_L1}.
In particular,
for $\mathring{\mathcal J}^m$, we first use the observation that $(t^{**}_1-s)/(t^*_m-s)\le t_1^{**}/t^*_m\simeq t_1/t_m$ for $s\in(0,t^{**}_1)$.
Then for  $\mathring{\mathcal J}^m$ and ${\mathcal J}^m$,
it is helpful to respectively use the substitutions $\hat s=s/t^{**}_1$ and $\hat s=s/t^*_m$,
while for ${\mathcal J}^m$ we also employ $(t^{**}_1/t^*_m)^{-\alpha/r}\sim (t_1/t_m)^{-\alpha/r}\sim (\tau_m/t_m)^{-\alpha}$ (also in view of \eqref{t_grid_gen}).

Combining the above observations  with \eqref{r_m_simplestA}  yields
$$
|r^m|\lesssim \max_{j=1,\ldots,m}\bigl\{\underbrace{(\tau_j/t_j)^{\fbox{\scriptsize$3$}-\alpha}}_{\simeq(\tau/t_j)^{(3-\alpha)/r}} (t_j/t_m)^{\alpha+1}
\,\psi^j\bigr\},
$$
where we also used $\tau_j/t_j\simeq(\tau/t_j)^{1/r}$ (in view of \eqref{t_grid_gen}).
The desired bound \eqref{tr_er_boundA} follows as $\tau\le t_j\le t_m$.
\end{proof}
\smallskip

{\it Proof of Theorem~\ref{lem_simplest_starA}.}
%\begin{proof}
Consider the error $e^m:=u(t_m)-U^m$, for which \eqref{simplestA} implies
$e^0=0$ and $\delta_{t}^{\alpha,*} e^m=r^m$  $\forall\,m\ge 1$, where
the truncation error $r^m$ is from Lemma~\ref{lem_truncA} and hence satisfies \eqref{tr_er_boundA}.
Furthermore, under the conditions \eqref{t_grid_gen} on the temporal mesh,
one has $\psi^1\lesssim 1$ (in view of ${\rm osc}\bigl(u, [0,t_2]\bigr)\le\int_0^{t_2}|\pt_s u|\,ds\lesssim t_2^{\alpha}$)
and $\psi^j\lesssim 1$ for $j\ge 2$ (in view of $s\simeq t_j$ for $s\in(t_{j-1},t_j)$ for this case).
Consequently, we arrive at
$$
|r^m|\lesssim (\tau/t_m)^{\gamma+1}\quad\forall\,m\ge 1,
\qquad\mbox{where}\;\;
\gamma+1:=\min\{\alpha+1,(3-\alpha)/r\}.
$$
The remainder of the proof employs Theorem~\ref{theo_barrierA} and closely follows the proof of Theorem~\ref{lem_simplest_star}.
In particular, the three cases $1\le r<3-\alpha$, $r=3-\alpha$ and $r>3-\alpha$
are considered separately, while $\tau\sim M^{-r}$ now implies $\tau^{(3-\alpha)/r}\sim M^{\alpha-3}$.
%
%\end{proof}
\endproof

\subsection{Error analysis for the Alikhanov-type semidiscretization in time}\label{ssec_A_semi}

Consider the semidiscretization of our problem~\eqref{problem} in time using $\delta_{t}^{\alpha,*}$ from \eqref{delta_t_def_ab}:%
%$\delta_t^\alpha$ from \eqref{delta_t_def}:
\begin{subequations}\label{semediscr_methodA_ab}
\beq\label{semediscr_methodA}
\delta_t^{\alpha,*} U^m +\LL U^{m,*}= f(\cdot,t^*_m)\;\;\mbox{in}\;\Omega,\quad U^m=0\;\;\mbox{on}\;\pt\Omega\quad\forall\,m=1,\ldots,M;\quad U^0=u_0.
\eeq
where, in view of \eqref{t_star}, we use a second-order discretization for $\LL u(\cdot, t_m^*)$ with
\beq\textstyle
U^{m,*}:=\frac12\alpha\,U^{m-1}+\bigl(1-\frac12\alpha\bigr)U^m.
\eeq
\end{subequations}

To simplify the presentation, here we shall consider only standard graded temporal meshes, which clearly satisfy both the condition from Remark~\ref{rem_JSC} and
 \eqref{t_grid_gen}. We shall also make some simplifying assumptions of $\LL$.
\smallskip

\begin{lemma}[Stability for fractional parabolic case]\label{lem_semi_stab}
Given $\gamma\in\R$,
let $\{t_j=T(j/M)^r\}_{j=0}^M$
for some $1\le r\le (3-\alpha)/\alpha$
if $\gamma>\alpha-1$ or for some $r\ge1$ if $\gamma\le \alpha-1$.
Also, let $\LL$ of \eqref{LL_def} be independent of $t$ with $b_k=0$ $\forall\,k$.
Then for $\{U^j\}_{j=0}^M$ from \eqref{semediscr_methodA_ab} one has
\beq\label{semi_stab}
\left.\begin{array}{c}
\|f(\cdot,t^*_j)\|_{L_2(\Omega)}
 \lesssim (\tau/ t_j)^{\gamma+1}
\\[0.2cm]
\forall j\ge1,\quad U^0=0\;\;\mbox{in}\;\bar\Omega
\end{array}\right\}
\;\;\Rightarrow\;\;
 \|U^j\|_{L_2(\Omega)}\lesssim
{\mathcal V}^j(\tau;\gamma)\quad\forall\, j\ge 1,
\eeq
where ${\mathcal V}^j={\mathcal V}^j(\tau;\gamma)$ is defined in \eqref{main_stab}.
\end{lemma}
\smallskip

\begin{proof}
(i)
Throughout the proof, we shall use the notation
$$
\delta_t^{\alpha,*} U^m=\sum_{j=0}^m\kappa_{m,j}U^j, \qquad f^j:=f(\cdot,t^*_j),\qquad \rho_j:=\tau_{j+1}/\tau_j,
$$
where, in view of Remark~\ref{rem_JSC}, $\kappa_{m,m}>0$, while $\kappa_{m,j}\le0$ $\forall\,j<m$.
An inspection of some arguments in \cite{ChenMS_JSC} shows (see Remark~\ref{rem_chenMSt_parabolic} below for further details)
that there exists a constant $c_0=c_0(\alpha,r)\in(0,1)$ such that
$\frac12\alpha\kappa_{m,m}<c_0(1-\frac12\alpha)|\kappa_{m,m-1}|$ \mbox{$\forall\, m\ge2$}.
%$$
%{\textstyle\frac12}\alpha\le c_0\bigl(1-{\textstyle\frac12}\alpha\bigr){|\kappa_{m,m-1}|}/{\kappa_{m,m}}.
%$$
Next, we claim that there is a sufficiently large $1\le K\lesssim 1$ (where $K=K(\alpha,r)$ is independent of $M$) such that
\beq\label{A_K_choice}
{\textstyle\frac12}\alpha\le \kappa_{m,m}^{-1/2}\,\kappa^{-1/2}_{m-1,m-1}\, \bigl(1-{\textstyle\frac12}\alpha\bigr){|\kappa_{m,m-1}|}\qquad\forall m\ge K+1.
\eeq
Indeed, it suffices to check that $c^2_0\le \kappa_{m,m}/\kappa_{m-1,m-1}$, while, in view of $\rho_{j}\le\rho_{j-1}$ $\forall\,j\ge 2$  a calculation shows that
$\tau_{m-1}^{\alpha}\kappa_{m-1,m-1}\le \tau_{m}^{\alpha}\kappa_{m,m}$, hence it suffices to check that $c_0^{2/\alpha}\le \tau_{m-1}/\tau_{m}$,
which can be ensured by choosing $K=K(\alpha, r)$ sufficiently large
(see the proof of \cite[Corollary~3.3]{NK_L2} for further details).

We shall consider the cases $K=1$ and $K>1$ separately in parts (ii) and (iii).%
\smallskip

(ii) Suppose $K=1$ in \eqref{A_K_choice}.
Then
\beq\label{w_prob}
\delta_t^{\alpha,*} w^m\le \|f^m\|_{L_2(\Omega)},\quad
w^m:=\sqrt{\displaystyle\|U^m\|_{L_2(\Omega)}^2+\kappa_{m,m}^{-1}\bigl(1-{\textstyle\frac12}\alpha\bigr)\langle\LL U^m, U^m\rangle}\,.
\eeq
Indeed, in view of \eqref{A_K_choice}, taking
the inner product of the equation from \eqref{semediscr_methodA} with
 $U^m$, one gets
$$
\kappa_{m,m} (w^m)^2\le |\kappa_{m,m-1}| w^m w^{m-1}+\sum_{j=1}^{m-2}|\kappa_{m,j}|\underbrace{\langle U^{m},U^j\rangle}_{\le w^m w^j}
+\underbrace{\langle U^m ,f^m\rangle}_{\le w^m\|f^m\|_{L_2(\Omega)}}\!\!\!\!.
$$
Dividing by $w^m$, we arrive at \eqref{w_prob}.
Now an application of Theorem~\ref{theo_barrierA}(iii) yields
$w^j\lesssim{\mathcal V}^j$, and hence
the desired result.%
\smallskip

(iii) Suppose that $1<K\lesssim 1$. We  imitate part (ii) in the proof of \cite[Theorem~3.2]{NK_L2}.
First, for $m\le K+1$, using $\tau_m\simeq \tau_1$ and \eqref{semediscr_methodA_ab}, one gets
$%\|U^m\|_{L_2(\Omega)}+\tau_1^\alpha \langle\LL U^m, U^m\rangle
w^m\lesssim w^{m-1}+\sum_{j=0}^{m-2}\|U^j\|_{L_2(\Omega)}+\tau_1^\alpha\|f^m\|_{L_2(\Omega)}$.
Here
$\|f^m\|_{L_2(\Omega)}\lesssim 1$, so
$\|U^m\|_{L_2(\Omega)}%+\tau_1^\alpha \langle\LL U^m, U^m\rangle
\lesssim \tau_1^\alpha\simeq {\mathcal V}^m$ $\forall\,m\le K+1$.
{\color{blue}Now \eqref{semediscr_methodA_ab} also implies that $\|\LL U^m\|_{L_2(\Omega)}\lesssim 1$ for $m\le K+1$.}

It remains to estimate the values of $\{\mathring{U}^j\}_{j=0}^M:=\{0,\ldots,0,U^{K+1},\ldots,U^M\}$
(i.e. $\mathring{U}^j$ is set to $0$ for $j\le K$ and to $U^j$ otherwise).
Note that $\delta_t^{\alpha,*} \mathring{U}^m+\LL \mathring{U}^{m,*}=0$ for $m\le K$
and ${\color{blue}\|\delta_t^{\alpha,*} \mathring{U}^{K+1}+\LL \mathring{U}^{K+1,*}\|_{L_2(\Omega)}}\lesssim 1$. Consider $m\ge K+2$. Then, by \eqref{delta_t_def}, one has
$\delta_t^{\alpha,*} \mathring{U}^m=\delta_t^{\alpha,*} {U}^m-D^\alpha_t \Pi^m[U-\mathring{U}](t^*_m)$.
As $\Pi^m[U-\mathring{U}]$ has support on $(0,t_{K+1})$, vanishes at $0$ and at $t_{K+1}\le t_m^*$, while its absolute value remains $\lesssim \tau_1^\alpha$,
so, recalling \eqref{CaputoEquiv} and applying an integration by parts yields
${\color{blue}\|{\color{black}D^\alpha_t \Pi^m[U-\mathring{U}](t_m)}\|_{L_2(\Omega)}}\lesssim \tau_1^\alpha\int_0^{t_{K+1}}(t^*_m-s)^{-\alpha-1}ds\lesssim (\tau_1/t_m)^{\alpha+1}$
(where we also used $t_{K+1}\simeq \tau_1$ and $t_m^*-s\gtrsim t_m^*-t_{K+1}\simeq t_m$).
Consequently, one concludes that 
$\color{blue}\|\delta_t^{\alpha,*} \mathring{U}^m+\LL \mathring{U}^{m,*}\|_{L_2(\Omega)}$, for $m\ge K+2$, is $\lesssim (\tau_1/t_m)^{\gamma+1}$ if $\gamma\le\alpha$
and $\lesssim (\tau_1/t_m)^{\alpha+1}$ otherwise.
%The latter result is also valid for $m=K+1$, which is shown  using a similar, but somewhat more intricate argument. Now we replace $\Pi^m[U-\mathring{U}](t_m)$ (as it does not vanish at $t_m^*$)
%by $I_2[U-\mathring{U}]$, where $I_2$ is from the proof of Lemma~\ref{lem_stability1_A} and is applied to $U-\mathring{U}$ set to $0$ at $t_m^*$.
%Then,
%in view of \eqref{nu_m},
%$\delta_t^{\alpha,*} \mathring{U}^m=\delta_t^{\alpha,*} {U}^m-D^\alpha_t I_2[U-\mathring{U}](t^*_m)$, and
%the desired bound on $|\delta_t^{\alpha,*} \mathring{U}^{K+1}|$ follows.

Finally, we restrict the problem for $\{\mathring{U}^j\}_{j=K-1}^M$ to the mesh $\{t_j\}_{j=K-1}^M$ and
note that for the Alikhanov-type operator $\mathring{\delta}_t^{\alpha,*} $ associated with the latter mesh one gets
$\mathring{\delta}_t^{\alpha,*}\mathring{U}^K=0 $ and $\mathring{\delta}_t^{\alpha,*}\mathring{U}^m={\delta}_t^{\alpha,*}\mathring{U}^m $ for $m\ge K+1$.
Now, in view of \eqref{A_K_choice}, an application of the result of part (ii) yields $ \|\mathring{U}^j\|_{L_2(\Omega)}\lesssim
{\mathcal V}^j$, which leads to the desired bound on $\|\mathring{U}^j\|_{L_2(\Omega)}$.
\end{proof}
\smallskip

\begin{remark}\label{rem_chenMSt_parabolic}
Comparing our notation $\{\kappa_{m,j}\}$ with \cite[(7)]{ChenMS_JSC}, the relation \cite[(15)]{ChenMS_JSC}
can be rewritten as
$\sigma|\kappa_{m,m-1}|-(1-\sigma)\kappa_{m,m}>0$, where $\sigma=1-\frac12\alpha$, while a sufficient condition for the latter is given by \cite[(17)]{ChenMS_JSC}
and is satisfied by our mesh.
Furthermore, an inspection of the proof of \cite[Lemma~4]{ChenMS_JSC} shows  such that
%under the same sufficient condition
%
%one in fact has $c_0\sigma\kappa_{m,m-1}-(1-c_0\sigma)\kappa_{m,m}>0$.
%Indeed,
in the second relation in \cite[(41)]{ChenMS_JSC}, one can include a constant factor $c_1(\alpha, \sigma\bar \rho)\in(0,1)$ in the right-hand side,
where $\bar\rho:=\max \rho_j=\rho_1$ on our  mesh.
(The latter observation can be shown by inspecting the proof of \cite[Lemma~2]{ChenMS_JSC} and replacing the the piecewise-constant approximation of $(t_{k+\sigma}-\eta)^{-\alpha}$ by a piecewise-linear one.)
Then (under the same sufficient condition)
one obtains a stronger version of \cite[(15)]{ChenMS_JSC}:
$\sigma^*|\kappa_{m,m-1}|-(1-\sigma^*)\kappa_{m,m}>0$ with $(2\sigma^*-1)/\sigma^*:=c_1(2-\sigma)/\sigma$, i.e. $\sigma^*=\sigma^*(\alpha, \sigma \bar\rho)<\sigma$.
Consequently, $(1-\sigma)\kappa_{m,m}<c_0\sigma|\kappa_{m,m-1}|$, where $c_0:=\frac{\sigma^*}{1-\sigma^*}\frac{1-\sigma}{\sigma }\in(0,1)$.
Recalling that $\sigma=1-\frac12\alpha$,
we conclude that $\forall\,m\ge2$ one has
$\frac12\alpha\kappa_{m,m}<c_0(1-\frac12\alpha)|\kappa_{m,m-1}|$.
\end{remark}
\smallskip

\begin{theorem}\label{theo_semidiscrA}
Let $\{t_j=T(j/M)^r\}_{j=0}^M$
for some $r\ge1$.
%Let the temporal mesh satisfy
% \eqref{t_grid_gen}
%with $r\ge 1$.
Suppose that $u$ is from \eqref{problem},\eqref{LL_def}, where
 $\LL$ of \eqref{LL_def} is independent of $t$ with $b_k=0$ $\forall\,k$.
Also, let $\|\partial_t^l u (\cdot, t)\|_{L_2(\Omega)}\lesssim 1+t^{\alpha-l}$ for $l = 1,3$
and $\|\partial_t^2\LL u (\cdot, t)\|_{L_2(\Omega)}\lesssim 1+t^{\alpha-2}$
$\forall\,t\in(0,T]$.
Then for  $\{U^m\}$ from \eqref{semediscr_methodA_ab}, one has
$$%\beq%\label{error}
\|u(\cdot,t_m)-U^m\|_{L_2(\Omega)}\lesssim {\mathcal E}^{m,**}:= {\mathcal E}^{m,*}+
M^{-2}\!\left\{\!\!\begin{array}{cl}
t_m^{2\alpha-2/r}&\mbox{if}\;\;2/r<\alpha+1\\
0&\mbox{otherwise}
\end{array}\right.
\quad\forall\,m\ge1,
$$%\eeq
where ${\mathcal E}^{m,*}$ is from \eqref{E_cal_mA}.
\end{theorem}
\smallskip

\begin{proof}
For the error $e^m:= u(\cdot,t_m)-U^m$, using \eqref{problem} and \eqref{semediscr_methodA_ab}, one immediately gets $e^0=0$ and $\forall\,m\ge1$
\beq\label{error_A_par}
\delta_t^{\alpha,*} e^m +\LL e^{m,*}=r^m+R^m, \qquad\mbox{where}\quad
\textstyle R^m:= \LL u^{m,*}-\LL u(\cdot,t^*_m)
\eeq
with the notation $u^{m,*}:=\frac12\alpha\,u(\cdot,t_{m-1})+\bigl(1-\frac12\alpha\bigr)u(\cdot,t_m)$,
and the truncation error $r^m$ from Lemma~\ref{lem_truncA} that satisfies
\eqref{tr_er_boundA}.
So under our assumptions of $u$ one has $\|R^m\|_{L_2(\Omega)}\lesssim \tau_m^2 t_m^{\alpha-2}\simeq \tau^{2/r}t_m^{\alpha-2/r}$ in view of \eqref{t_grid_gen},
and also $\|r^m\|_{L_2(\Omega)}\lesssim(\tau/t_m)^{\gamma+1}$, where
$\gamma+1:=\min\{\alpha+1,(3-\alpha)/r\}$ (see the proof  of Theorem~\ref{lem_simplest_starA}).

If  $2/r\ge\alpha+1\ge \gamma+1$, then  $\|R^m\|_{L_2(\Omega)}\lesssim (\tau/t_m)^{2/r}\le  (\tau/t_m)^{\gamma+1}$,
so  Lemma~\ref{lem_semi_stab} yields
$\|e^m\|_{L_2(\Omega)}\lesssim {\mathcal V}^m(\tau;\gamma)$.
Otherwise, %i.e. if $r> 2/(\alpha+1)$,
Lemma~\ref{lem_semi_stab} yields
$\|e^m\|_{L_2(\Omega)}\lesssim {\mathcal V}^m(\tau;\gamma)+\tau^\alpha{\mathcal V}^m(\tau;\gamma')$,
where $\gamma':=2/r-\alpha-1<0$ implies $\tau^\alpha{\mathcal V}^m(\tau;\gamma')=\tau^{2/r}t_m^{2\alpha-2/r}$, where $\tau^{2/r}\simeq M^{-2}$.
Finally,  imitating the proof of Theorem~\ref{lem_simplest_starA}, one gets
${\mathcal V}^m(\tau;\gamma)\lesssim{\mathcal E}^{m,*}$, so combining our findings we arrive at the desired error bound.
\end{proof}
\smallskip

\begin{remark}[Convergence in positive time]\label{rem_positive_time_A_par}
Consider $t_m\gtrsim 1$. Then, in view of Remark~\ref{rem_positive_time_A},
${\mathcal E}^{m,**}\simeq M^{-r}$ for $r\le 2$. Otherwise, $2/r<1<\alpha+1$ so ${\mathcal E}^{m,**}\simeq M^{-2}$.
In summary, ${\mathcal E}^{m,**}\simeq M^{-\min \{ r,2\}}$, and
$r=2$ yields
the optimal convergence rate $2$.

By contrast, \cite[Theorem~3.9]{arxiv18_liao_et_al} (obtained by means of a discrete Gr\"{o}nwall inequality \cite{sinum19_liao_et_al}) gives a somewhat similar, but less sharp error bound for graded meshes, as (in our notation)
it involves the term $O(\tau^\alpha)=O(M^{-\alpha r})$, so, e.g., for $r=2$ it
gives a considerably less sharp convergence rate of only $2\alpha$.
%it  requires $r=2/\alpha$ to attain the optimal convergence rate in positive time.
For $r=1$, we have ${\mathcal E}^{m,**}\simeq M^{-1}$, so our error bound  is again sharper than those in \cite[Theorem~3.9]{arxiv18_liao_et_al}.%
\end{remark}%
\smallskip

\begin{remark}[Global convergence]\label{rem_global_time_A_par}
In view of Remark~\ref{rem_global_time_A},
 $%\displaystyle
\max_{m\ge1}{\mathcal E}^{m,*}\simeq M^{-\min \{\alpha r,3-\alpha\}}$.
If $\alpha r\ge1$, then $2/r\le 2\alpha<\alpha+1$ so
$\max_{m\ge1}{\mathcal E}^{m,**}\simeq\max_{m\ge1}{\mathcal E}^{m,*}+M^{-2}\simeq M^{-\min \{\alpha r,2\}}$.
Otherwise, $\alpha r<1$ implies $\max_{m\ge1}{\mathcal E}^{m,*}\simeq M^{-\alpha r}$, while
 $\max_{m\ge1} M^{-2}t_m^{2\alpha-2/r}$ is attained at $m=1$ and is $\simeq \tau^{2\alpha}\simeq M^{-2\alpha r}$,
 so $\max_{m\ge1}{\mathcal E}^{m,**}%\simeq\max_{m\ge1}{\mathcal E}^{m,*}
 \simeq M^{-\alpha r}$.
 %
% Combining our findings, $\max_{m\ge1}{\mathcal E}^{m,**}\simeq M^{-\min \{\alpha r,2\}}$ for any $r\ge1$.
%
%
Consequently, Theorem~\ref{theo_semidiscrA} yields the global error bound
$|u(t_m)-U^m|\lesssim M^{-\min \{\alpha r,2\}} $.
This implies that
the optimal grading parameter for global accuracy is $r=2/\alpha$.
Note that a similar global error bound was obtained in \cite{ChenMS_JSC}.
\end{remark}
%\smallskip

\subsection{Alikhanov-type finite element discretizations}\label{ssec_A_full}

Discretize \eqref{problem}--\eqref{LL_def}, posed in a general bounded  Lipschitz domain  $\Omega\subset\R^d$,
 by applying
 a standard Galerkin finite element spatial approximation, described in \S\ref{ssec_L1_full_FEM}, to the temporal semidiscretization~\eqref{semediscr_methodA_ab}.
%
%  A Lagrange finite element space $S_h \subset H_0^1(\Omega)\cap C(\bar\Omega)$ of fixed degree $\ell\ge 1$, % \Bbbk
% relative to a
% quasiuniform simplicial triangulation
% of $\Omega$, is employed, as in \cite[\S5]{NK_MC_L1}.
 Then for the full discretization one can easily generalize the stability result given by Lemma~\ref{lem_semi_stab}.
 Furthermore, for $e_h^m:=\RR_h u(t_m)-u_h^m\in S_h$, where
 $u_h^m\in S_h$ is the finite element solution at time $t_m$, and
 $\RR_h u(t)\in S_h$ is the  Ritz projection of $u(\cdot,t)$,
 a standard calculation (see, e.g., \cite[Theorem~5.1]{NK_MC_L1}) yields
 $$
 \langle \delta_t^{\alpha,*} e_h^m,v_h\rangle +\AA (e_h^{m,*},v_h) =\langle \delta_t^{\alpha,*}\!\rho+r^m+R^m ,v_h\rangle\quad \forall v_h\in S_h.
 $$
 Here $\AA(\cdot,\cdot)$ is the standard bilinear form associated with $\LL$, and $\rho(\cdot, t):=\RR_h u(t)-u(\cdot, t)$.
 So,
 under the conditions of Theorem~\ref{theo_semidiscrA}  and assuming that $u_h^0=\RR_h u(0)$,
 one gets
the following version of \cite[Theorem~5.5]{NK_L2}:
\begin{align*}%\label{err_A_fem}
\|u(\cdot,t_m)-u_h^m\|_{L_2(\Omega)}\lesssim{}&{\mathcal E}^{m,**}
+\|\rho(\cdot, t_m)\|_{L_2(\Omega)}+t_m^\alpha\sup_{t\in(0,t_{m})}\!\!\bigl\{t^{1-\alpha} \|\pt_t\rho(\cdot, t)\|_{L_2(\Omega)}\bigr\}
%\qquad\forall\,m\ge 1,
\end{align*}
$\forall\,m\ge 1$, where ${\mathcal E}^{m,**}$ is from  Theorem~\ref{theo_semidiscrA}.
%\\
%
%Here $u_h^m\in S_h$ is the finite element solution at time $t_m$, ${\mathcal E}^m$ is from \eqref{E_cal_m}, and
%$\rho(\cdot, t):=\RR_h u(t)-u(\cdot, t)$ is the error of
%the standard Ritz projection $\RR_h u(t)\in S_h$ of $u(\cdot,t)$.
Under additional realistic assumptions on $u$, the final two terms in the above error estimate can be bounded by $O( h^{\ell+1})$,
where $h$ is the triangulation diameter \cite[\S5]{NK_MC_L1}.

\section{Numerical results}\label{sec_Num}

\subsection{Fractional parabolic test with finite elements}\label{ssec_num1}
Our first fractional-order parabolic test problem  is \eqref{problem} with $\LL=-(\pt_{x_1}^2+\pt_{x_2}^2)$, posed
in the domain $\Omega\times[0,1]$ %(see Fig.\,\ref{fig_mesh}, left)
from \cite[\S7]{NK_MC_L1}
with $\pt\Omega$  parameterized by
$x_1(l):=\frac23R\cos\theta$ and $x_2(l):=R\sin\theta$, where
$R(l): =  0.8 + \cos^2\! l$ and
$\theta(l) := l + e^{(l-5)/2}\sin(l/2)\sin l$ for $l\in[0,2\pi]$. %; see \cite[\S7]{NK_MC_L1}.
We choose  $f$, as well as the initial and non-homogeneous boundary conditions, so that the unique exact solution
$u=t^{\alpha}[1+\ln(x-y/3+7)]$.
This problem is discretized in space (with an obvious modification for the case of non-homogeneous boundary conditions)
using
lumped-mass linear finite elements %(described in \S\ref{ssec_lumped})
on quasiuniform Delaunay triangulations of $\Omega$ (with DOF denoting the number of degrees of freedom in space).
{The errors will be computed in the approximate $L_2(\Omega)$ norm as
$\|u_h-u^I\|_{L_2(\Omega)}$, where $u^I\in S_h$ is the piecewise-linear interpolant in $\Omega$.}
%
%The graded temporal mesh $\{t_j=T(j/M)^r\}_{j=0}^M$ is used in all numerical experiments.
All numerical experiments will use the graded temporal mesh
$\{t_j=T(j/M)^r\}_{j=0}^M$.

For the L1 method,
we have the error bounds \eqref{L1_semi_error} and \eqref{L1_fem_error}.
%applied to fractional-order parabolic problems.
These error bounds are consistent with the
numerical rates of convergence %for this method,
given in \cite{gracia_etal_cmame} for errors in positive time and  $r=1$, as well as those in
\cite{stynes_etal_sinum17,NK_MC_L1} for
 errors in the maximum norm in time and
various $r$. % These results  are consistent with our error bounds.
Additionally, consider the case $r>2-\alpha$, for which our error bounds predict the optimal convergence rate of $2-\alpha$ with respect to time
at $t\gtrsim 1$ (see Remark~\ref{rem_positive_time}).
This agrees with the numerical convergence rates given in Table~\ref{t2} for the L1 method with $r=(2-\alpha)/0.9$.

For the Alikhanov method, %applied to fractional-order parabolic problems
we have the error bounds of Theorem~\ref{theo_semidiscrA} and \S\ref{ssec_A_full}.
Note that they are consistent with the
numerical rates of convergence given in~\cite{ChenMS_JSC} for
 errors in the maximum norm in time and
various $r$.
Additionally, here we numerically investigate the case $r=2$, for which our error bounds predict the optimal convergence rate $2$ with respect to time
at $t\gtrsim 1$ (see Remark~\ref{rem_positive_time_A_par}).
This clearly agrees  with the numerical convergence rates given in
Table~\ref{t2} for the the Alikhanov method.

%{\color{red}For these two methods, we also show the errors in the maximum $L_2(\Omega)$ norm for a large fixed $M$
%in
%Fig.\,\ref{fig_mesh} (right).
%As $M$ is quite large, the error bounds from \eqref{L1_fem_error} and  \S\ref{ssec_A_full} imply that
%the errors are expected to be $\lesssim h^2\simeq \mbox{DOF}^{-1}$, which is again consistent with our numerical results.
%}

%  \begin{figure}[t!]%[b!]
%% \vspace*{-0.5pc}
%\begin{center}
%\includegraphics[height=0.38\textwidth]{domain}\hfill%\includegraphics[height=0.38\textwidth]{L2_method_pic_}
%\end{center}
%\vspace{-0.3cm}
% \caption{\label{fig_mesh}\it\small
% Fractional-order parabolic test problem:
% Delaunay triangulation of $\Omega$ with DOF=172 (left),
%maximum $L_2(\Omega)$ and $L_\infty(\Omega)$ errors for $\alpha=0.5$, $r=(3-\alpha)/\alpha$ and $M=2048$.}
% \end{figure}

 {
\begin{table}[h!]
%\beforecaption={\tabcolsep=0pt}
\begin{center}
\caption{
%Alikhanov method applied to the fractional-order parabolic test problem:
% $L_2(\Omega)$ errors at $t=1$ (odd rows) and
%computational rates $q$ in $M^{-q}$ (even rows) for
% $r=2$ and spatial \color{red}DOF=255435%
Fractional-order parabolic test problem from \S\ref{ssec_num1}:
 $L_2(\Omega)$ errors at $t=1$ (odd rows) and
computational rates $q$ in $M^{-q}$ (even rows) for the
L1 method with $r=(2-\alpha)/.9$ and the
Alikhanov method with
 $r=2$, spatial DOF=398410
 }
\label{t2}%\raggedright
\vspace{-0.1cm}
\tabcolsep=4pt
{\small
\begin{tabular}{lrrrrrrrr}
\hline
\strut\rule{0pt}{9pt}
&\multicolumn{4}{l}{L1 method,~ $r=\frac{2-\alpha}{.9}$}&
\multicolumn{4}{l}{~~~Alikhanov method,~ $r=2$}\\[2pt]
\hline
\strut\rule{0pt}{9pt}&
$M=2^6$& $M=2^7$& $M=2^8$& $M=2^9$&{}~~~$M=2^6$& $M=2^7$& $M=2^8$& $M=2^9$\\%&$M=2048$\\
\hline
$\alpha=0.3\;\;$

%L1
%&8.346e-05	&2.579e-05	&7.984e-06	&2.481e-06	
%&1.694	&1.692	&1.686

%A:
%M_vec = 2.^(6:9) = 64   128   256   512
%N=500
%dof = 398410
%&5.790e-06	&1.251e-06	&2.828e-07	&7.036e-08	
%&2.210	&2.146	&2.007	

&8.35e-5	&2.58e-5	&7.98e-6	&2.48e-6
        &5.79e-6	&1.25e-6	&2.83e-7	&7.04e-8\\
&1.69	&1.69	&1.69
        &&2.21	&2.15	&2.01\\[3pt]	
$\alpha=0.5$
%
%L1
%/0.9

&2.66e-4	&9.47e-5	&3.38e-5	&1.20e-5	%&4.30e-06	
        &7.10e-6	&1.58e-6	&3.67e-7	&9.15e-8\\
&1.49	&1.49	&1.49	%&1.49	
        &&2.16	&2.11	&2.00\\[3pt]
%-------------
%&2.659e-04	&9.470e-05	&3.376e-05	&1.204e-05	&4.297e-06	
%&1.490	&1.488	&1.487	&1.486	
%A:
%M_vec = 2.^(6:9) = 64   128   256   512
%N=500
%dof = 398410
%&7.101e-06	&1.584e-06	&3.672e-07	&9.153e-08\\	
%&2.164	&2.109	&2.004\\[3pt]

%&7.10e-6	&1.58e-6	&3.67e-7	&9.15e-8\\
%&&2.16	&2.11	&2.00\\[3pt]		

%
$\alpha=0.7$
%	
%L1
%&5.606e-04	&2.299e-04	&9.442e-05	&3.879e-05	
%&1.286	&1.284	&1.284
%A:
%M_vec = 2.^(6:9) = 64   128   256   512
%N=500
%dof = 398410
%&7.383e-06	&1.719e-06	&4.089e-07	&1.026e-07	
%&2.103	&2.072	&1.995	
%
&5.61e-4	&2.30e-4	&9.44e-5	&3.88e-5
        &7.38e-6	&1.72e-6	&4.09e-7	&1.03e-7\\	
&1.29	&1.28	&1.28	
        &&2.10	&2.07	&1.99\\	
\hline
\end{tabular}}
\end{center}
\end{table}
}

\subsection{Fractional parabolic test with finite differences}\label{ssec_num_FD}
To test the error bound
\eqref{FD_error} given in \S\ref{ssec_FD} for finite difference discretizations in space combined with the L1 scheme in time,
we shall employ another test problem. Consider
\eqref{problem}
with $\LL=-(\pt_{x_1}^2+\pt_{x_2}^2)+(1+x_1+x_2+t)$,
 the initial condition $u_0=\sin x_1\,\sin x_2$,
 and $f=x_1(\pi-x_1)x_2(\pi-x_2)(1+t^4) + t^2$,
 posed
in the domain $\Omega\times[0,1]$ with the square spatial domain $\Omega=(0,\pi)^2$
(this test is a modification of \cite[Example 6.2]{stynes_etal_sinum17}). %
 The spatial mesh was a uniform tensor product mesh of size $h=\pi/N$ (i.e. with $N$ equal mesh intervals in each coordinate direction). As the exact solution is unknown, the errors were computed using  the two-mesh principle.

  We focus on the most interesting case of
  the graded temporal mesh
 with $r>2-\alpha$, for which our error bound~\eqref{FD_error} predicts the optimal convergence rate of $2-\alpha$ with respect to time
at $t\gtrsim 1$ (in view of Remark~\ref{rem_positive_time}).
This clearly agrees with the numerical convergence rates given in Table~\ref{t2plus}  for
the grading parameter
$r=(2-\alpha)/0.9$.
 %The errors, in the nodal maximum norm, and the corresponding convergence rates are given in Table~\ref{t2plus}.
 %We observe that the convergence rates clearly agree with the error bound~\eqref{FD_error}.

%\end{example}
%\newpage

\begin{table}[t!]
\begin{center}
\caption{
{Fractional-order parabolic test problem from \S\ref{ssec_num_FD}:
maximum  nodal errors at $t=1$ (odd rows) and
computational rates $q$ in $M^{-q}$ or $N^{-q}$ (even rows) for the L1 method with $r=(2-\alpha)/.9$}
}
\label{t2plus}
\vspace{-0.1cm}
\tabcolsep=4pt
{\small
\begin{tabular}{lrrrrrrrr}
\hline
\strut\rule{0pt}{9pt}
&\multicolumn{4}{l}{errors and convergence rates in time}&
\multicolumn{4}{l}{~~~errors and convergence rates in space}\\[0pt]
\strut\rule{0pt}{9pt}
&\multicolumn{4}{l}{$N = M$}&
\multicolumn{4}{l}{~~~$M = N^2%{\lceil 2/(2-\alpha) \rceil}
$}\\[2pt]
\hline
\strut\rule{0pt}{9pt}&
$M=2^5$& $M=2^6$& $M=2^7$& $M=2^8$&{}~~~$N=2^3$&~$N=2^4$&~$N=2^5$&~$N=2^6$\\
\hline
$\alpha=0.3\;\;$

&6.99e-4	&2.30e-4	&7.45e-5	&2.39e-5
       &2.81e-3	&7.36e-4	&1.87e-4	&4.86e-5\\
&1.60	&1.63	&1.64
        &&1.93	&1.97	&1.95\\[3pt]

$\alpha=0.5$

&1.54e-3	&5.75e-4	&2.10e-4	&7.59e-5		
        &2.87e-3	&7.34e-4	&1.84e-4	&4.86e-5\\
&1.43	&1.45	&1.47	
        &&1.97	&1.99	&1.92\\[3pt]

$\alpha=0.7$

&3.05e-3	&1.28e-3	&5.29e-4	&2.17e-4
        &3.15e-3	&7.84e-4	&1.91e-4	&4.86e-5\\	
&1.25	&1.27	&1.28	
        &&2.01	&2.04	&1.97\\

%M=N^2/4
%&3.70e-03	&8.32e-04	&1.97e-04	&4.85e-05	
%&5.24e-03	&1.03e-03	&2.19e-04	&5.00e-05	
%&8.08e-03	&1.60e-03	&3.14e-04	&6.45e-05	
%&2.15	&2.08	&2.02	
%&2.35	&2.23	&2.13	
%&2.34	&2.35	&2.29	

%M=N^2/2
%&3.00e-3	&7.58e-4	&1.90e-4	&4.86e-5	
%&3.40e-3	&8.03e-4	&1.93e-4	&4.86e-5	
%&4.52e-3	&9.81e-4	&2.20e-4	&5.10e-5	
%&1.98	&2.00	&1.96	
%&2.08	&2.06	&1.99	
%&2.20	&2.15	&2.11	

%M=N^2/1

\hline
\end{tabular}}
\end{center}
\end{table}

\subsection{L1 method: pointwise sharpness of the error estimate for the  initial-value problem}\label{ssec_num_L1}
Here, to demonstrate the sharpness of the error estimate \eqref{E_cal_m} given by Theorem~\ref{lem_simplest_star} for the L1 method,
we consider the simplest initial-value fractional-derivative test problem \eqref{simplest} %without spatial derivatives
with the simplest typical exact solution $u(t):=t^\alpha$.
Table~\ref{t0_positive_time} shows the errors and the corresponding convergence rates at $t=1$, which agree with~\eqref{E_cal_m}, in view of
Remark~\ref{rem_positive_time}.
In particular, the latter implies that the errors are
$\lesssim M^{-\min \{ r,2-\alpha\}}$ for $r\neq 2-\alpha$ and $\lesssim M^{-(2-\alpha)}\ln M$ for $r= 2-\alpha$.
The maximum errors and corresponding convergence rates for various $\alpha$ and $r$ are given
in  \cite{stynes_etal_sinum17,NK_MC_L1}, and they
%in Table~\ref{t0_global}
 confirm the conclusions of Remark~\ref{rem_global_time},
which predicts from the pointwise  bound~\eqref{E_cal_m} that the global errors are
$\lesssim M^{-\min \{\alpha r,2-\alpha\}}$.

\begin{figure}[h!]%[b!]
% \vspace*{-0.5pc}
\begin{center}
\includegraphics[height=0.30\textwidth]{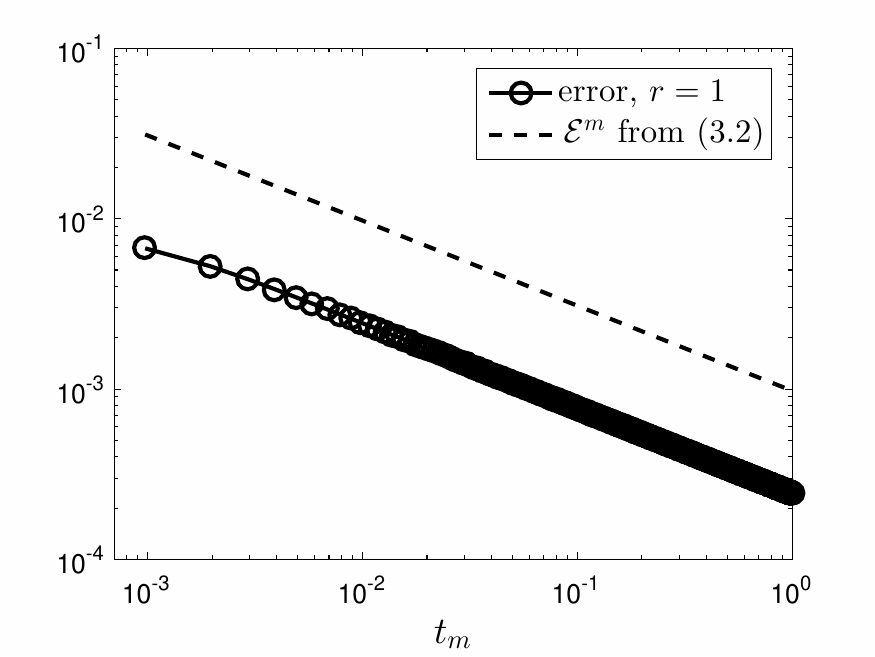}~~\includegraphics[height=0.30\textwidth]{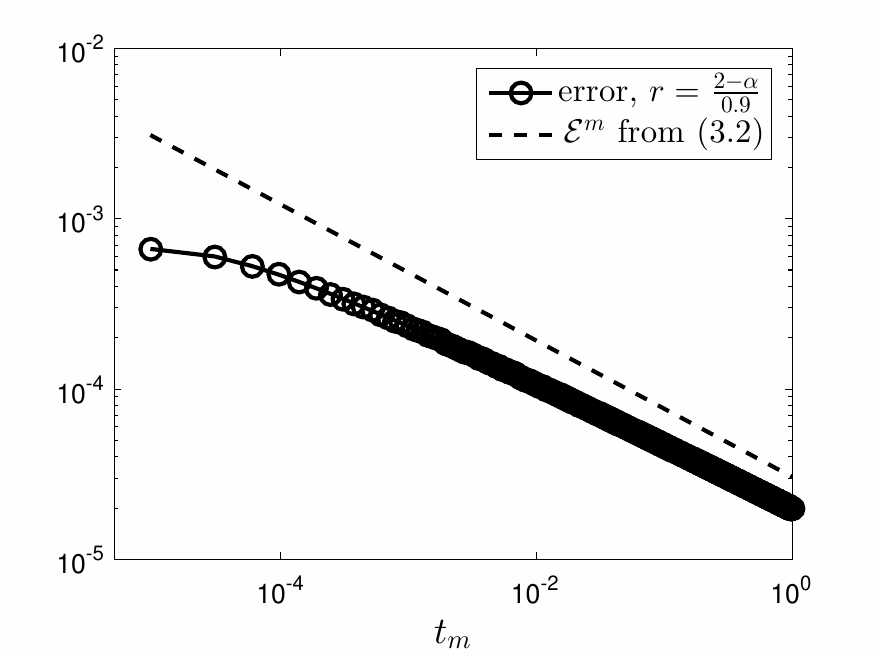}%
\\%
\includegraphics[height=0.30\textwidth]{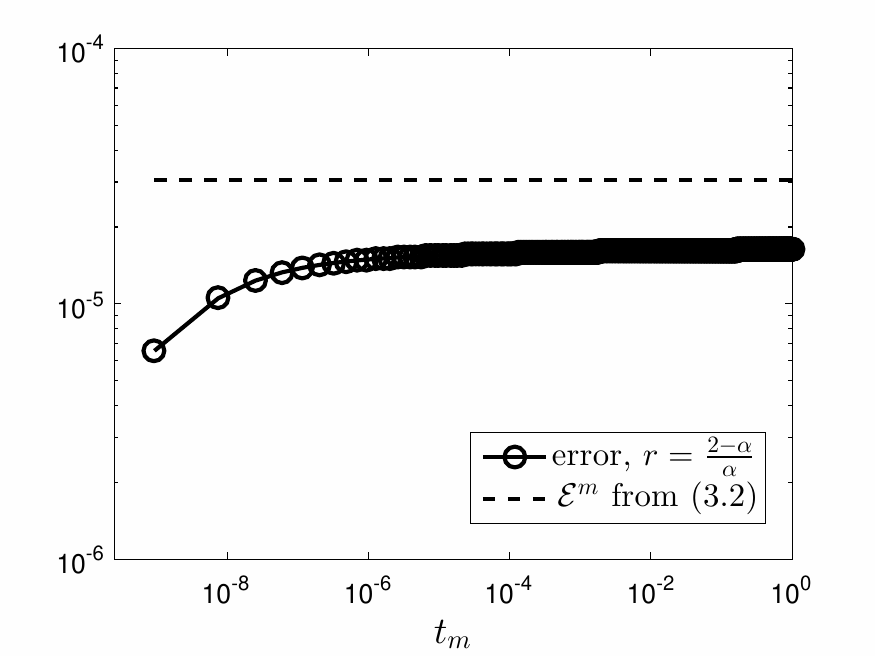}~~\includegraphics[height=0.30\textwidth]{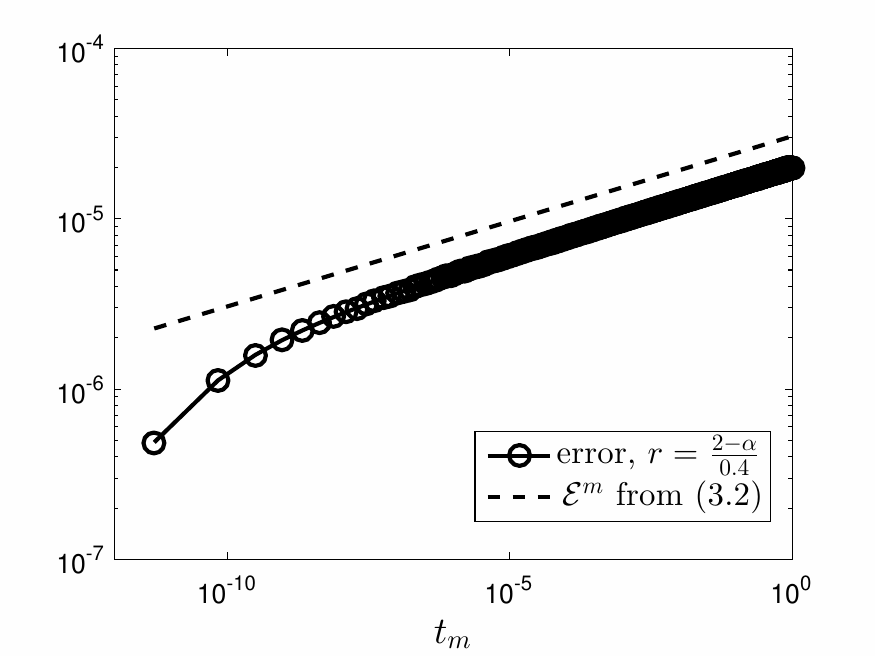}%
\end{center}
\vspace{-0.3cm}
 \caption{\label{pointwise_fig}\it\small
L1 method applied to the initial-value test problem: pointwise errors for $\alpha=0.5$ and $M=1024$, cases $r=1$, $r=(2-\alpha)/0.9$, $r=(2-\alpha)/\alpha$ and $r=(2-\alpha)/0.4$.}
 \end{figure}

Furthermore, in Fig.~\ref{pointwise_fig}, the pointwise errors for various $r$ are compared with the pointwise theoretical error bound~\eqref{E_cal_m},
and again, with the exception of a few initial mesh nodes,  we observe remarkably good agreement.
Note that Fig.~\ref{pointwise_fig} only addresses the case $\alpha=0.5$, but for other values of $\alpha$ we observed similar consistency of \eqref{E_cal_m}
with the actual pointwise errors.

{
\begin{table}[t!]
%\beforecaption={\tabcolsep=0pt}
\begin{center}
\caption{%
%\parbox[t]{9cm}{
L1 method applied to the initial-value test problem:
  errors at $t=1$ (odd rows) and
computational rates $q$ in $M^{-q}$ (even rows) for
 $r=1$, $r=2-\alpha$ and
 $r=(2-\alpha)/.95$}%}
\label{t0_positive_time}%\raggedright
\tabcolsep=5pt
\vspace{-0.1cm}
{\small
\begin{tabular}{lrrrrrrr}
\hline
\strut\rule{0pt}{9pt}&&
$\;\;\;M=2^7$& $M=2^9$&$M=2^{11}$&
 $M=2^{13}$&$M=2^{15}$&$M=2^{17}$\\
\hline
$r=1$
&$\alpha=0.3$
&1.182e-3	&2.939e-4	&7.333e-5	&1.832e-5	&4.578e-6	&1.144e-6\\	
&&1.004	&1.001	&1.001	&1.000	&1.000\\[2pt]
&$\alpha=0.5$
&1.953e-3	&4.883e-4	&1.221e-4	&3.052e-5	&7.629e-6	&1.907e-6\\
&&1.000	&1.000	&1.000	&1.000	&1.000\\[2pt]
&$\alpha=0.7$	
&2.489e-3	&6.433e-4	&1.642e-4	&4.163e-5	&1.050e-5	&2.640e-6\\	
&&0.976	&0.985	&0.990	&0.994	&0.996\\
%
   %-------------------------------
\hline\strut\rule{0pt}{9pt}
$r=2-\alpha\;$
&$\alpha=0.3$
&1.201e-4	&1.310e-5	&1.401e-6	&1.477e-7	&1.540e-8	&1.592e-9\\	
&&1.598	&1.612	&1.623	&1.631	&1.637\\[2pt]	
&$\alpha=0.5$
&5.039e-4	&7.407e-5	&1.063e-5	&1.500e-6	&2.089e-7	&2.878e-8\\
&&1.383	&1.400	&1.413	&1.422	&1.430\\[2pt]
&$\alpha=0.7$
&1.267e-3	&2.495e-4	&4.782e-5	&8.986e-6	&1.663e-6	&3.042e-7\\	
&&1.172	&1.192	&1.206	&1.217	&1.225\\
%
 %-------------------------------
\hline
\strut\rule{0pt}{9pt}
$r=\frac{2-\alpha}{.95}$
&$\alpha=0.3$
&1.035e-4	&1.074e-5	&1.094e-6	&1.098e-7	&1.092e-8	&1.076e-9\\	
&&1.634	&1.648	&1.658	&1.665	&1.671\\[2pt]
&$\alpha=0.5$
&4.469e-4	&6.276e-5	&8.609e-6	&1.161e-6	&1.546e-7	&2.039e-8\\	
&&1.416	&1.433	&1.445	&1.454	&1.461\\[2pt]
&$\alpha=0.7$
&1.143e-3	&2.164e-4	&3.984e-5	&7.192e-6	&1.279e-6	&2.250e-7\\	
&&1.201	&1.221	&1.235	&1.245	&1.254\\
\hline
\end{tabular}}
\end{center}
\end{table}
}

{
\begin{table}[t!]
%\beforecaption={\tabcolsep=0pt}
\begin{center}
\caption{
%\parbox[t]{9cm}{
Alikhanov method applied to the initial-value test problem:
  errors at $t=1$ (odd rows) and
computational rates $q$ in $M^{-q}$ (even rows) for
 $r=1$, $r=2$ and
 $r=(3-\alpha)/.95$}%}
\label{t0_positive_time_A}%\raggedright
\tabcolsep=5pt
\vspace{-0.1cm}
{\small
\begin{tabular}{lrrrrrrr}
\hline
\strut\rule{0pt}{9pt}&&
 %$\;\;\;M=2^5$&
$\;\;\;M=2^6$& $M=2^8$&$M=2^{10}$&
 $M=2^{12}$&$M=2^{14}$&$M=2^{16}$\\
\hline
$r=1$
&$\alpha=0.3$
&1.325e-3	&3.306e-4	&8.260e-5	&2.065e-5	&5.162e-6	&1.290e-6\\[2pt]	
&&1.002	&1.000	&1.000	&1.000	&1.000\\
&$\alpha=0.5$
&1.530e-3	&3.819e-4	&9.543e-5	&2.386e-5	&5.964e-6	&1.491e-6\\[2pt]	
&&1.001	&1.000	&1.000	&1.000	&1.000\\
&$\alpha=0.7$
&1.236e-3	&3.087e-4	&7.715e-5	&1.929e-5	&4.821e-6	&1.205e-6\\	
&&1.001	&1.000	&1.000	&1.000	&1.000\\	
%
   %-------------------------------
\hline\strut\rule{0pt}{9pt}
$r=2$
&$\alpha=0.3$
&3.891e-5	&2.446e-6	&1.530e-7	&9.560e-9	&5.975e-10	&3.734e-11\\
&&1.996	&1.999	&2.000	&2.000	&2.000\\[2pt]
&$\alpha=0.5$
&6.079e-5	&3.940e-6	&2.502e-7	&1.576e-8	&9.885e-10	&6.190e-11\\	
&&1.974	&1.988	&1.995	&1.997	&1.999\\[2pt]	
&$\alpha=0.7$
&6.450e-5	&4.436e-6	&2.936e-7	&1.902e-8	&1.216e-09	&7.720e-11\\	
&&1.931	&1.959	&1.974	&1.984	&1.989\\  	
%
 %-------------------------------
\hline
\strut\rule{0pt}{9pt}
$r=\frac{3-\alpha}{.95}$
&$\alpha=0.3$
&1.085e-5	&3.241e-7	&8.953e-9	&2.363e-10	&6.058e-12	&1.509e-13\\
&&2.532	&2.589	&2.622	&2.643	&2.664\\[2pt]
&$\alpha=0.5$
&2.710e-5	&1.057e-6	&3.839e-8	&1.337e-9	&4.529e-11	&1.517e-12\\	
&&2.340	&2.392	&2.422	&2.442	&2.450\\[2pt]	
&$\alpha=0.7$
&3.962e-5	&2.017e-6	&9.638e-8	&4.431e-9	&1.986e-10	&8.791e-12\\	
&&2.148	&2.194	&2.221	&2.240	&2.249\\
\hline
\end{tabular}}
\end{center}
\end{table}
}

{
\begin{table}[t!]
%\beforecaption={\tabcolsep=0pt}
\begin{center}
\caption{
%\parbox[t]{9cm}{
Alikhanov method applied to the initial-value test problem:
 maximum nodal errors (odd rows) and
computational rates $q$ in $M^{-q}$ (even rows) for
 $r=1$, $r=2/\alpha$ and $r=(3-\alpha)/\alpha$}%}
\label{t0_global_A}%\raggedright
\tabcolsep=5pt
\vspace{-0.1cm}
{\small
\begin{tabular}{lrrrrrrr}
\hline
\strut\rule{0pt}{9pt}&&
%$M=2^5$&
$\;\;\;M=2^6$& $M=2^8$&$M=2^{10}$&
 $M=2^{12}$&$M=2^{14}$&$M=2^{16}$\\
\hline
$r=1$
&$\alpha=0.3$
&2.477e-2	&1.634e-2	&1.078e-2	&7.115e-3	&4.694e-3	&3.097e-3\\	
&&0.300	&0.300	&0.300	&0.300	&0.300\\[2pt]
&$\alpha=0.5$	
&1.164e-2	&5.819e-3	&2.909e-3	&1.455e-3	&7.273e-4	&3.637e-4\\
&&0.500	&0.500	&0.500	&0.500	&0.500\\[2pt]		
&$\alpha=0.7$
&3.919e-3	&1.485e-3	&5.627e-4	&2.132e-4	&8.079e-5	&3.061e-5\\
&&0.700	&0.700	&0.700	&0.700	&0.700\\	
%
  %-------------------------------
\hline\strut\rule{0pt}{9pt}
$r=\frac{2}{\alpha}$
&$\alpha=0.3$
&5.865e-5	&3.665e-6	&2.291e-7	&1.432e-8	&8.949e-10	&5.593e-11\\
&&2.000	&2.000	&2.000	&2.000	&2.000\\[2pt]
&$\alpha=0.5$
&5.250e-5	&3.281e-6	&2.051e-7	&1.282e-8	&8.011e-10	&5.007e-11\\
&&2.000	&2.000	&2.000	&2.000	&2.000\\[2pt]
&$\alpha=0.7$
&4.232e-5	&2.645e-6	&1.653e-7	&1.033e-8	&6.458e-10	&4.036e-11\\	
&&2.000	&2.000	&2.000	&2.000	&2.000\\	
%
 %-------------------------------
\hline
\strut\rule{0pt}{9pt}
$r=\frac{3-\alpha}{\alpha}$
&$\alpha=0.3$
&5.505e-5	&1.659e-6	&4.472e-8	&1.142e-9	&2.833e-11	&6.923e-13	\\
&&2.526	&2.607	&2.646	&2.667	&2.677\\[2pt]
&$\alpha=0.5$
&3.976e-5	&1.379e-6	&4.508e-8	&1.439e-9	&4.542e-11	&1.425e-12	\\
&&2.425	&2.467	&2.485	&2.493	&2.497\\[2pt]
&$\alpha=0.7$
&3.425e-5	&1.498e-6	&6.307e-8	&2.619e-9	&1.083e-10	&4.469e-12	\\
&&2.257	&2.285	&2.295	&2.298	&2.299\\
\hline
\end{tabular}}
\end{center}
\end{table}
}

\subsection{Alikhanov method: pointwise sharpness of the error estimate for the  initial-value problem}\label{ssec_num_A}
Next, we turn to the Alikhanov method and,
to demonstrate the sharpness of the error estimate (\ref{E_cal_mA}) given by Theorem~\ref{lem_simplest_starA}, %for the Alikhanov method,
%we
consider the simplest initial-value fractional-derivative test problem \eqref{simplestA} %without spatial derivatives
with the same simplest typical exact solution $u(t):=t^\alpha$.
Table~\ref{t0_positive_time_A} shows the errors and the corresponding convergence rates at $t=1$, which agree with~(\ref{E_cal_mA}), in view of
Remark~\ref{rem_positive_time_A}.
In particular, the latter implies that the errors are
$\lesssim M^{-\min \{ r,3-\alpha\}}$ for $r\neq 3-\alpha$.
The maximum errors and corresponding convergence rates given in Table~\ref{t0_global_A}
clearly confirm the conclusions of Remark~\ref{rem_global_time_A},
which predicts from the pointwise  bound~(\ref{E_cal_mA}) that the global errors are
$\lesssim M^{-\min \{\alpha r,3-\alpha\}}$.

Note that, similarly to Fig.~\ref{pointwise_fig}, we observed the pointwise behaviour of the errors consistent with (\ref{E_cal_mA});
see also \cite[Fig.~6.2]{NK_L2} for similar graphs of pointwise errors of an L2-type method.

\section*{Acknowledgents}
The authors are grateful to Prof. Martin Stynes for his helpful comments which inspired the %inclusion in the manuscript of Section~4 on
extension of our analysis to
the Alikhanov method.

%\newpage

\end{document}